\documentclass[10pt]{article}
\usepackage{color}
\usepackage{amssymb,amsmath,graphicx,acronym,amsfonts}
\newtheorem{corollary}{Corollary}[section]
\newtheorem{theorem}{Theorem}[section]
\newtheorem{lemma}{Lemma}[section]
\newtheorem{definition}{Definition}[section]
\newtheorem{proposition}{Proposition}[section]
\newtheorem{example}{Example}[section]
\newtheorem{assum}{Assumption}[section]
\newtheorem{algo}{Algorithm}[section]
\newtheorem{Remark}{Remark}[section]

\def\bc{\begin{corl}}
\def\bc{\end{corl}}
\def\ba{\begin{algo}}
\def\ea{\end{algo}}
\def\br{\begin{Remark}}
\def\er{\end{Remark}}
\def\bs{\begin{assum}}
\def\es{\end{assum}}
\def\bt{\begin{theorem}}
\def\et{\end{theorem}\vskip 3pt}
\def\bl{\begin{lemma}}
\def\el{\end{lemma}}
\def\ep{\end{proposition}}
\def\bp{\begin{proposition}}
\def\qed{\hfill{$\Box$}\vskip 5pt}
\def\be{\begin{example}}
\def\ee{\end{example}}
\def\bd{\begin{definition}}
\def\ed{\end{definition}}
\def\bc{\begin{corollary}}
\def\ec{\end{corollary}}
\def\proof{\noindent\it Proof. \hspace{1mm}\rm}
\input epsf
\begin{document}
\title{\bf \Large SOS Tensor Decomposition: Theory and Applications}
\author{Haibin Chen\thanks{Department of Applied Mathematics, The Hong Kong Polytechnic University, Hung Hom,
Kowloon, Hong Kong. Email: chenhaibin508@163.com. This author's work was supported by the Natural Science Foundation of China (11171180).}, \quad
Guoyin Li\thanks{Department of Applied Mathematics, University of New South Wales, Sydney 2052, Australia. E-mail: g.li@unsw.edu.au
(G. Li). This author's work was partially supported by Australian Research Council.},\quad
Liqun Qi
\thanks{Department of Applied Mathematics, The Hong Kong Polytechnic University, Hung Hom,
Kowloon, Hong Kong. Email: maqilq@polyu.edu.hk. This author's work
was supported by the Hong Kong Research Grant Council (Grant No.
PolyU 502111, 501212, 501913 and 15302114).} }

\date{}
\maketitle

\begin{abstract}
In this paper, we examine structured tensors which have sum-of-squares (SOS) tensor decomposition, and study the SOS-rank of SOS tensor decomposition.
We first show that several classes of even order symmetric structured tensors available in the literature have SOS tensor decomposition.
These include positive Cauchy tensors, weakly diagonally dominated tensors,
$B_0$-tensors,  double $B$-tensors, quasi-double $B_0$-tensors, $MB_0$-tensors, $H$-tensors, absolute tensors
of positive semi-definite $Z$-tensors and extended $Z$-tensors.
We also examine the SOS-rank of SOS tensor decomposition and the SOS-width for SOS tensor cones. The SOS-rank provides the minimal number of squares in the SOS tensor
decomposition, and, for a given SOS tensor cone, its SOS-width is the maximum possible SOS-rank for all the tensors in this cone. We first
deduce an upper bound for general tensors that have SOS decomposition and the SOS-width for general SOS tensor cone using the known results in the literature of polynomial theory.
Then, we provide an explicit sharper
estimate for the SOS-rank of SOS tensor decomposition with bounded exponent and identify the SOS-width for the tensor cone consisting of all tensors with bounded exponent that have SOS decompositions.
Finally, as applications, we show how the SOS tensor decomposition can be used to compute the minimum $H$-eigenvalue of an even order symmetric extended $Z$-tensor and test the positive definiteness of an associated multivariate form. Numerical experiments are also provided to show the efficiency of the proposed numerical methods ranging from small size to large size numerical examples.
\medskip

\noindent{\bf Keywords:} structured tensor, SOS tensor decomposition, positive semi-definite tensor, SOS-rank, $H$-eigenvalue.
\vskip 6pt

\noindent{\bf AMS Subject Classification(2000):} 90C30,  15A06.

\end{abstract}

\newpage
\section{Introduction}
Tensor decomposition is an important research area, and it has found numerous applications in data mining \cite{Kolda2006,Kolda2008,Kolda2009}, computational neuroscience \cite{Field1991,Comon2014}, and statistical learning for latent variable models \cite{Anan14}.
An important class of tensor decomposition is sum-of-squares (SOS) tensor decomposition. It is known that to determine a given even order symmetric tensor is positive
semi-definite or not is an NP-hard problem in general.  On the other hand, an interesting feature of SOS tensor decomposition is checking whether a given even order symmetric tensor has SOS decomposition or not can be verified by solving a semi-definite programming problem (see for example \cite{Hu14}), and hence, can be validated efficiently. SOS tensor decomposition has a close connection with SOS polynomials, and SOS polynomials are very
important in polynomial theory \cite{Ch77,Ch95,Habi39,H88,Power98,Rez00} and polynomial optimization \cite{KLYZ,JB01,M09,Le14,parrilo03,shor98}.
It is known that an even order symmetric tensor having SOS decomposition
 is positive semi-definite, but the converse is not true in general. Recently, a few classes of structured tensors such as B tensors \cite{qi14} and diagonally dominated tensor \cite{Qi05},
have been shown to be positive semi-definite in the even order symmetric case. It then raises a natural and interesting question:  Will these structured tensors admit an SOS decomposition? Providing an answer for this question is important because this will enrich the theory of SOS tensor decomposition, achieve
a better understanding for these structured tensors, and lead to efficient numerical methods for solving problems involving these structured tensors.

In this paper, we make the following contributions in answering the above theoretical question and providing applications on important numerical problems involving structured tensors:
\begin{itemize}
\item[{\rm (1)}] We first show that several classes of symmetric structured tensors available in the literature have SOS decomposition when the order is even. These classes include positive Cauchy tensors, weakly diagonally dominated tensors,
$B_0$-tensors, double $B$-tensors, quasi-double $B_0$-tensors, $MB_0$-tensors, $H$-tensors, absolute tensors
of positive semi-definite $Z$-tensors and extended $Z$-tensors.
\item[{\rm (2)}] Secondly, we examine the SOS-rank for tensors with SOS decomposition and the SOS-width for SOS tensor cones.
The SOS-rank of tensor $\mathcal{A}$ is defined to be the minimal number of the squares which appear in the sums-of-squares decomposition
of the associated homogeneous polynomial of $\mathcal{A}$, and, for a given SOS tensor cone, its SOS-width is the maximum possible SOS-rank for all the tensors in this cone.  We deduce an upper bound for the SOS-rank of general SOS tensor decomposition and the SOS-width for the general SOS tensor cone using the known result in
polynomial theory \cite{Ch95}. We then provide a sharper explicit upper bound of the SOS-rank for tensors with bounded exponent and identify the exact SOS-width for the  cone consists of all such tensors with SOS decomposition.
\item[{\rm (3)}] Finally, as applications, we show how the derived SOS tensor decomposition can be used to compute the minimum $H$-eigenvalue of an even order symmetric extended $Z$-tensor and test the positive definiteness of an associated multivariate form. Numerical experiments are also provided to show the efficiency of the proposed numerical method ranging from small size to large size numerical examples.
\end{itemize}


The rest of this paper is organized as follows. In Section 2, we recall some basic definitions and facts for tensors and polynomials.
We also present some properties of SOS tensor cone and its duality. In Section 3, we present SOS decomposition property for various classes of structured tensors. In Section 4, we study the SOS-rank of SOS tensor decomposition and SOS-width for a given SOS tensor cone.
In particular, we examine SOS tensor decomposition with bounded exponents and the SOS-width
of the cone constituted by all such tensors with SOS decomposition, and provide their sharper explicit estimate.
In Section 5, as applications for the derived SOS decomposition of structure tensors, we show that
the minimum $H$-eigenvalue of an even order extended $Z$-tensor can be computed via polynomial optimization technique.
Accordingly, this also leads to an efficient test for positive definiteness of an associated multivariate form.
Numerical experiments are also provided to illustrate the significance of the result. Final remarks and some questions are listed in Section 6.

Before we move on, we briefly mention the notation that will be used in the sequel. Let $\mathbb{R}^n$
be the $n$ dimensional real Euclidean space and the set consisting of all positive integers
is denoted by $\mathbb{N}$. Suppose $m, n\in \mathbb{N}$ are two natural numbers. Denote $[n]=\{1,2,\cdots,n\}$. Vectors are denoted by bold lowercase letters i.e. ${\bf x},~ {\bf y},\cdots$, matrices are denoted by capital letters i.e. $A, B, \cdots$, and tensors are written as calligraphic capitals such as
$\mathcal{A}, \mathcal{T}, \cdots.$ The $i$-th unit coordinate vector in $\mathbb{R}^n$ is denoted by ${\bf e_i}$.
If the symbol $|\cdot|$ is used on a tensor $\mathcal{A}=(a_{i_1 \cdots i_m})_{1\leq i_j\leq n}$, $j=1,\cdots,m$, it denotes another
tensor $|\mathcal{A}|=(|a_{i_1 \cdots i_m}|)_{1\leq i_j\leq n}$, $j=1,\cdots,m$.

\setcounter{equation}{0}
\section{Preliminaries}
A real $m$th order n-dimensional tensor $\mathcal{A}=(a_{i_1i_2\cdots i_m})$ is a multi-array
of real entries $a_{i_1i_2\cdots i_m}$, where $i_j \in [n]$ for $j\in [m]$.
If the entries $a_{i_1i_2\cdots i_m}$ are invariant under any permutation of
their indices, then tensor $\mathcal{A}$ is called a symmetric tensor. In this paper, we always consider
symmetric tensors defined in $\mathbb{R}^n$.  The identity tensor $\mathcal{I}$ with order $m$ and dimension $n$ is given by $\mathcal{I}_{i_1\cdots i_m}=1$ if $i_1=\cdots=i_m$ and $\mathcal{I}_{i_1\cdots i_m}=0$ otherwise.

We first fix some symbols and recall some basic facts of tensors and polynomials.  Let $m,n \in \mathbb{N}$. Consider
$S_{m,n}:=\{\mathcal{A}: \mathcal{A} \mbox{ is an } m\mbox{th-order }  n\mbox{-dimensional} \mbox{  symmetric tensor}\}.$
Clearly, $S_{m,n}$ is a vector space under the addition and multiplication
defined as below: for any $t \in \mathbb{R}$,
$\mathcal{A}=(a_{i_1 \cdots i_m})_{1 \le i_1,\cdots,i_m
\le n}$ and $\mathcal{B}=(b_{i_1 \cdots i_m})_{1 \le
i_1,\cdots,i_m \le n},$
\[
\mathcal{A}+\mathcal{B}=(a_{i_1 \cdots i_m}+b_{i_1 \cdots i_m})_{1
\le i_1,\cdots,i_m \le n} \mbox{ and } t
\mathcal{A}=(ta_{i_1 \cdots i_m})_{1 \le i_1,\cdots,i_m
\le n}.
\]
For each $\mathcal{A}, \mathcal{B} \in S_{m,n}$, we define the inner
product by
\[
\langle
\mathcal{A},\mathcal{B}\rangle:=\sum_{i_1,\cdots,i_m=1}^{n}a_{i_1 \cdots
i_m}b_{i_1 \cdots i_m}.
\]
The corresponding norm is defined by $\displaystyle
\|\mathcal{A}\|=\left(\langle
\mathcal{A},\mathcal{A}\rangle\right)^{1/2}=\left(\sum_{i_1,\cdots,i_m=1}^{n}(a_{i_1
\cdots i_m})^2\right)^{1/2}$.
For a vector ${\bf x}\in
\mathbb{R}^n$, we use $x_i$ to denote its $i$th component.
Moreover, for a  vector ${\bf x}\in
\mathbb{R}^n$, we use ${\bf x}^m$ to denote the $m$th-order $n$-dimensional symmetric rank one tensor induced by ${\bf x}$, i.e.,
\[
({\bf x}^m)_{i_1 i_2\cdots i_m}=x_{i_1}x_{i_2}\cdots x_{i_m}, \ \forall \, i_1,\cdots,i_m \in \{1,\cdots,n\}.
\]

We note that an $m$th order $n$-dimensional symmetric tensor uniquely defines an $m$th degree homogeneous
polynomial $f_{\mathcal{A}}$ on $\mathbb{R}^n$: for all ${\bf x}=(x_1,\cdots,x_n)^T
\in \mathbb{R}^n$,
\begin{equation}\label{e21}
f_{\mathcal{A}}({\bf x})= \mathcal{A}{\bf x}^m=\sum_{i_1,i_2,\cdots, i_m\in [n]}a_{i_1i_2\cdots i_m}x_{i_1}x_{i_2}\cdots x_{i_m}.
\end{equation}
Conversely, any $m$th degree homogeneous
polynomial function $f$ on $\mathbb{R}^n$ also uniquely corresponds a symmetric tensor. Furthermore,
a tensor $\mathcal{A}$ is called positive semi-definite (positive definite) if $f_{\mathcal{A}}({\bf x}) \geq 0$ ($f_{\mathcal{A}}({\bf x})> 0$)
for all ${\bf x}\in \mathbb{R}^n$ (${\bf x}\in \mathbb{R}^n \backslash \{\bf 0\}$).

{We now recall the following definitions on eigenvalues and eigenvectors for a tensor \cite{Lim05, Qi05}.
\bd\label{def21}  Let $\mathbb{C}$ be the complex field. Let $\mathcal{A}=(a_{i_1i_2\cdots i_m})$ be an order $m$ dimension $n$ tensor. A pair $(\lambda, {\bf x})\in \mathbb{C}\times \mathbb{C}^n\setminus \{0\}$ is called an
eigenvalue-eigenvector pair of tensor $\mathcal{A}$, if they satisfy	
$$
\mathcal{A}{\bf x}^{m-1}=\lambda {\bf x}^{[m-1]},
$$
where $\mathcal{A}{\bf x}^{m-1}$ and ${\bf x}^{[m-1]}$ are all n dimensional column vectors given by
$$\mathcal{A}{\bf x}^{m-1}=\left(\sum_{i_2,\cdots,i_m=1}^n a_{ii_2\cdots i_m}x_{i_2}\cdots x_{i_m} \right)_{1\leq i\leq n}$$ and ${\bf x}^{[m-1]}=(x_1^{m-1},\ldots,x_n^{m-1})^T \in \mathbb{R}^n$.
\ed
If the eigenvalue $\lambda$ and the eigenvector ${\bf x}$ are real, then
$\lambda$ is called an $H$-eigenvalue of $\mathcal{A}$ and ${\bf x}$ is its corresponding $H$-eigenvector \cite{Qi05}.
An important fact which will be used frequently later on is that
an even order symmetric tensor is positive semi-definite (definite) if and only if all $H$-eigenvalues of the tensor are nonnegative (positive).}

Suppose that $m$ is even. In (\ref{e21}),
if $f_{\mathcal{A}}({\bf x})$ is a sums-of-squares (SOS) polynomial, then we say $\mathcal{A}$ has an {\bf SOS tensor decomposition} (or an SOS decomposition, for simplicity). It is clear that a tensor with SOS decomposition and an SOS polynomial must have even degree.
If a given tensor has SOS decomposition, then
the tensor is positive semi-definite, but not vice versa. Next, we recall a useful lemma which provides a test for verifying whether a homogeneous
polynomial is a sums-of-squares polynomial or not. To do this, we introduce some basic notions.

For all ${\bf x}\in \mathbb{R}^n$, consider a homogeneous polynomial $f({\bf x})=\sum_{\alpha}f_{\alpha}{\bf x}^{\alpha}$ with degree $m$ ($m$ is an even number), where $\alpha=(\alpha_1,\cdots,\alpha_n) \in (\mathbb{N} \cup \{0\})^n$, ${\bf x}^{\alpha}=x_1^{\alpha_1}\cdots x_n^{\alpha_n}$ and $|\alpha|:=\sum_{i=1}^n \alpha_i=m$.
Let $f_{m,i}$ be the coefficient associated with $x_i^{m}$. Let ${\bf e_i}$ be the $i$th unit vector and let
\begin{equation}\label{e22} \Omega_f=\{\alpha=(\alpha_1,\cdots,\alpha_n) \in (\mathbb{N} \cup \{0\})^n: f_{\alpha} \neq 0 \mbox{ and } \alpha \neq m \, {\bf e_i}, \ i=1,\cdots,n\}.\end{equation}
Then, $f$ can be decomposed as
$f({\bf x})=\sum_{i=1}^n f_{m,i} x_i^{m}+\sum_{\alpha \in \Omega_f}f_{\alpha}{\bf x}^{\alpha}$. Recall that $2\mathbb{N}$ denotes the set consisting of all the even numbers. Define
$$
\hat{f}({\bf x})=\sum_{i=1}^n f_{m,i} x_i^{m}-\sum_{\alpha \in \Delta_f}|f_{\alpha}|{\bf x}^{\alpha},
$$
where
 \begin{equation}\label{e23}
\Delta_f:=\{\alpha=(\alpha_1,\cdots,\alpha_n) \in \Omega_f: f_{\alpha} < 0 \mbox{ or } \alpha \notin (2\mathbb{N} \cup \{0\})^n\}.
\end{equation}

\begin{lemma}\cite[Corollary 2.8]{FK11}\label{lema21}
Let $f$ be a homogeneous polynomial of degree
$m$, where $m$ is an even number. If $\hat{f}$ is a polynomial which always takes nonnegative values, then $f$ is a sums-of-squares polynomial.
\end{lemma}

\subsection*{SOS tensor cone and its dual cone}
In this part, we study the cone consisting of all tensors that have SOS decomposition, and its dual cone \cite{Luo}.  We use ${\rm SOS}_{m,n}$ to denote the cone consisting of all order $m$ and dimension $n$ tensors, which have SOS decomposition. The following simple lemma from \cite{Hu14} gives some basic properties of ${\rm SOS}_{m,n}$.
\begin{lemma}\label{lema31} {\rm (cf. \cite{Hu14})}
Let $m,n \in \mathbb{N}$ and $m$ be an even number. Then, ${\rm SOS}_{m,n}$ is a closed convex cone with dimension at most $I(m,n)=\binom {n+m-1} {m}.$
\end{lemma}

For a closed convex cone $C$, we  recall that the dual cone of $C$ in
$S_{m,n}$ is denoted by $C^{\oplus}$ and defined by
$C^{\oplus}=\{\mathcal{A} \in S_{m,n}: \langle
\mathcal{A},\mathcal{C}\rangle \ge 0$ for all $\mathcal{C} \in C\}$.
Let $\mathcal{M}=(m_{i_1,i_2,\cdots,i_m}) \in S_{m,n}$. We also define the symmetric tensor
${\rm sym}(\mathcal{M} \otimes \mathcal{M}) \in S_{2m,n}$ by $${\rm
sym}(\mathcal{M} \otimes \mathcal{M}){\bf x}^{2m}=
(\mathcal{M}{\bf x}^ m)^2=\sum_{1 \le i_1,\cdots,  i_m,
j_1,\cdots,j_m \le n}m_{i_1, \cdots ,i_m}m_{j_1,
\cdots ,j_m}x_{i_1}\cdots x_{i_m}x_{j_1}\cdots x_{j_m}.$$
Moreover, in the case where the degree $m=2$,
${\rm SOS}_{2,n}$ and its dual cone are equal, and both reduce to the cone of positive semidefinite $(n \times n)$ matrices. Therefore, to avoid triviality,
we consider the duality of the SOS tensor cone ${\rm SOS}_{m,n}$ in the case where $m$ is an even number with $m \ge 4$.

\begin{proposition}{\bf (Duality between tensor cones)}
Let $n \in \mathbb{N}$ and $m$ be an even number with $m \ge 4$. Then, we have
${\rm SOS}_{m,n}^{\oplus} = \{\mathcal{A} \in S_{m,n}: \langle \mathcal{A}, {\rm sym}(\mathcal{M} \otimes \mathcal{M})\rangle \ge 0, \,~ \forall~\mathcal{M} \in S_{\frac{m}{2},n}\}$ and ${\rm SOS}_{m,n} \nsubseteqq {\rm SOS}_{m,n}^{\oplus}$. 
\end{proposition}

\proof We define ${\rm SOS}_{m,n}^h$ to be the cone consisting of all $m$th-order $n$-dimensional symmetric tensors such that
 $f_{\mathcal{A}}({\bf x}):=\langle \mathcal{A}, {\bf x}^ m \rangle$ is a polynomial which can be written as sums of finitely many homogeneous polynomials. We now see that indeed ${\rm SOS}_{m,n}^{h} = {\rm SOS}_{m,n}$.  Clearly, ${\rm SOS}_{m,n}^{h} \subseteq {\rm SOS}_{m,n}$. To see the reverse inclusion, we let $\mathcal{A} \in {\rm SOS}_{m,n}$. Then, there exists $l \in \mathbb{N}$ and $f_1,\cdots,f_l$ are real polynomials with degree at most $\frac{m}{2}$ such that
 $\langle \mathcal{A}, {\bf x}^m \rangle=\sum_{i=1}^l f_i({\bf x})^2$. In particular, for all $t \ge 0$, we have
 \[
 t^m \, \langle \mathcal{A}, {\bf x}^m\rangle =\langle \mathcal{A}, (t{\bf x})^m \rangle=\sum_{i=1}^l f_i(t{\bf x})^2
 \]
 Dividing $t^m$ on both sides and letting $t \rightarrow +\infty$, we see that
 $\langle \mathcal{A}, {\bf x}^m\rangle =\sum_{i=1}^l f_{i,\frac{m}{2}} ({\bf x})^2$,
 where $f_{i,\frac{m}{2}}$ is the $\frac{m}{2}$th-power term of $f_i$, $i=1,\cdots,l$. This shows that $\mathcal{A} \in {\rm SOS}_{m,n}^{h}$. Thus, we have ${\rm SOS}_{m,n}^{h} = {\rm SOS}_{m,n}$.  It then follows that
\begin{eqnarray*}
\big({\rm SOS}_{m,n} \big)^{\oplus}  =  \big({\rm SOS}_{m,n}^h \big)^{\oplus}
&=&\{\mathcal{A} \in S_{m,n}: \langle \mathcal{A}, \mathcal{C}\rangle \ge 0 \mbox{ for all } \mathcal{C} \in {\rm SOS}^h_{m,n}\}\\
& = & \{\mathcal{A} \in S_{m,n}: \langle \mathcal{A}, \mathcal{C}\rangle \ge 0 \mbox{ for all } \mathcal{C}=\sum_{i=1}^l {\rm sym}(\mathcal{M}_i \otimes \mathcal{M}_i), \\
& & \ \ \ \ \ \ \ \ \ \ \ \ \ \ \ \ \ \ \mathcal{M}_i \in S_{\frac{m}{2},n}, i=1,\cdots,l\} \\
& = & \{\mathcal{A} \in S_{m,n}: \langle \mathcal{A}, {\rm sym}(\mathcal{M} \otimes \mathcal{M})\rangle \ge 0 \mbox{ for all }  \mathcal{M} \in S_{\frac{m}{2},n}\}.
\end{eqnarray*}
We now show that ${\rm SOS}_{m,n} \nsubseteqq {\rm
SOS}_{m,n}^{\oplus}$ if $m \ge 4$. Let
$f({\bf x})=x_1^4+x_2^4+\frac{1}{4}x_3^4+6x_1^2x_2^2+6x_1^2x_3^2+6x_2^2x_3^2$
and let $\mathcal{A} \in S_{4,3}$ be such that
$\mathcal{A}{\bf x}^4=f({\bf x})$. Then, $\mathcal{A}$ has an SOS
decomposition and $\mathcal{A}_{1,1,1,1}=\mathcal{A}_{2,2,2,2}=1$,
$\mathcal{A}_{3,3,3,3}=\frac{1}{4}$,
$\mathcal{A}_{1,1,3,3}=\mathcal{A}_{1,1,2,2}=\mathcal{A}_{2,2,3,3}=1$.
We now see that $\mathcal{A} \notin {\rm SOS}_{m,n}^{\oplus}$. To
see this, we only need to find $M \in S_{2,3}$ such that
$\langle \mathcal{A},{\rm sym}(M \otimes M)\rangle<0.$
To see this, let $M={\rm diag}(1,1,-4)$. Then, ${\rm sym}(M \otimes M){\bf x}^4=({\bf x}^TM{\bf x})^2=(x_1^2+x_2^2-4x_3^2)^2$. Direct verification shows that
${\rm sym}(M \otimes M){\bf x}^4=x_1^4+x_2^4+16x_3^4+2x_1^2x_2^2-8x_1^2x_3^2-8x_2^2x_3^2.$
So, ${\rm sym}(M \otimes M)_{1,1,1,1}={\rm sym}(M \otimes M)_{2,2,2,2}=1$, ${\rm sym}(M \otimes M)_{3,3,3,3}=16$, ${\rm sym}(M \otimes M)_{1,1,2,2}=\frac{1}{3}$, ${\rm sym}(M \otimes M)_{1,1,3,3}={\rm sym}(M \otimes M)_{2,2,3,3}=-\frac{4}{3}$. Therefore,
\[
\langle \mathcal{A},{\rm sym}(M \otimes M)\rangle=1+1+\frac{1}{4}\cdot 16+ 6\left(1 \cdot \frac{1}{3}\right)+6 \left(1 \cdot \left(-\frac{4}{3}\right)\right) +6 \left(1 \cdot \left(-\frac{4}{3}\right)\right)=-8 <0,
\]
and the desired results hold.\qed

{\bf Question:} 
 It is known from polynomial optimization (see \cite[ Proposition 4.9]{M09} or \cite{JB01}) that the dual cone of the cone consisting of all sums-of-squares polynomials (possibly nonhomogeneous) is the moment cone (that is, all the sequence whose associated moment matrix is positive semi-definite). Can we link the dual cone of ${\rm SOS}_{m,n}$ to the moment matrix? Can the membership problem of ${\rm SOS}_{m,n}^{\oplus}$ be solvable in polynomial time?

\setcounter{equation}{0}
\section{SOS Decomposition of Several Classes of Structured Tensors}
In this section, we examine the SOS decomposition of several classes of symmetric even order structured tensors,
such as positive Cauchy tensor, weakly diagonally dominated tensors, $B_0$-tensors, double $B$-tensors, quasi-double $B_0$-tensors, $MB_0$-tensors, $H$-tensors, absolute tensors of positive semi-definite $Z$-tensors and extended $Z$-tensors.

\subsection{Characterizing SOS decomposition for even order Cauchy tensors}
Symmetric Cauchy tensors was first studied in \cite{chen14}. Some checkable sufficient and necessary conditions for
an even order symmetric Cauchy tensor to be positive semi-definite or positive definite were provided in \cite{chen14}, which extends the matrix cases established in \cite{Fied10}.

Let ${\bf c}=(c_1,c_2,\cdots,c_n)^T\in \mathbb{R}^n$ with $c_{i_1}+c_{i_2}+\cdots+c_{i_m} \neq 0$ for all $i_j \in \{1,\ldots,n\}$, $j=1,\ldots,m$. Let the real tensor $\mathcal{C}=(c_{i_1i_2\cdots i_m})$ be defined by
$$
c_{i_1i_2\cdots i_m}=\frac{1}{c_{i_1}+c_{i_2}+\cdots+c_{i_m}},\quad j\in [m],~i_j \in [n].
$$
Then, we say that $\mathcal{C}$ is a symmetric {\bf Cauchy tensor} with order $m$ and dimension $n$ or simply a Cauchy tensor. The corresponding vector ${\bf c}\in \mathbb{R}^n$ is called the generating vector of $\mathcal{C}$.

To establish the SOS decomposition of Cauchy tensors, we will also need another class of tensors called completely positive tensor, which has an SOS tensor
decomposition in the even order case.

Tensor $\mathcal{A}$ is called a {\bf completely decomposable tensor} if there are
vectors ${\bf x}_j\in \mathbb{R}^n$,~$j\in [r]$ such that $\mathcal{A}$ can be written as sums of rank-one tensors generated by the vector ${\bf x}_j$, that is,
$$\mathcal{A}=\sum\limits_{j\in [r]}{\bf x}_j^m.$$
If ${\bf x}_j\in \mathbb{R}^n_+$ for all $j\in [r]$, then $\mathcal{A}$ is called a {\bf completely positive tensor} \cite{QXX}. It was shown that a strongly
symmetric, hierarchically dominated nonnegative tensor is a completely positive tensor \cite{QXX}.

We now characterize the SOS decomposition and completely positivity for even order Cauchy tensors.
\bt\label{them36}
Let ${\bf c}=(c_1,c_2,\cdots,c_n)^T\in \mathbb{R}^n$ with $c_{i_1}+c_{i_2}+\cdots+c_{i_m} \neq 0$ for all $i_j \in \{1,\ldots,n\}$, $j=1,\ldots,m$. Let $\mathcal{C}$ be a Cauchy tensor generated by ${\bf c}$ with even order $m$ and dimension $n$.
 Then, the following statements are equivalent.
 \begin{itemize}
\item[{\rm (i)}] the  Cauchy tensor $\mathcal{C}$ has an SOS tensor decomposition;
\item[{\rm (ii)}] the  Cauchy tensor $\mathcal{C}$ is positive semi-definite;
\item[{\rm (iii)}]  the generating vector of the Cauchy tensor $c_i$, $i \in [n]$, are all positive;
\item[{\rm (iv)}] the  Cauchy tensor $\mathcal{C}$ is a completely positive tensor.
\end{itemize}
\et

\proof Since $m$ is even, by definitions of completely positive tensor, SOS tensor decomposition and positive
semi-definite tensor, we can easily obtain ${\rm (i)}\Rightarrow{\rm (ii)}$ and ${\rm (iv)}\Rightarrow{\rm (i)}$. By Theorem 2.1 of \cite{chen14}, we know that $\mathcal{C}$ is positive semi-definite if
and only if $c_i>0, i\in [n]$, and hence, ${\rm (ii)} \Leftrightarrow {\rm (iii)}$ holds.
So, we only need to prove ${\rm (iii)}\Rightarrow{\rm (iv)}$, that is, any Cauchy tensors with positive generating vector is completely positive.

Assume ${\rm (iii)}$ holds.  Then, for any ${\bf x}\in \mathbb{R}^n$,
$$
\begin{aligned}
\mathcal{C} {\bf x}^m  =& \sum_{i_1,i_2,\cdots,i_m=1}^n \frac{x_{i_1}x_{i_2} \cdots x_{i_m}}{c_{i_1}+c_{i_2}+\cdots+c_{i_m}}  \\
= & \sum_{i_1,i_2,\cdots,i_m=1}^n  \left(\int_0^1  t^{c_{i_1}+c_{i_2}+\cdots+c_{i_m}-1}x_{i_1}x_{i_2} \cdots x_{i_m} dt \right) \\
=&\int_0^1  \left(\sum_{i_1,i_2,\cdots,i_m=1}^n t^{c_{i_1}+c_{i_2}+\cdots+c_{i_m}-1} x_{i_1}x_{i_2} \cdots x_{i_m} \right)dt  \\
=& \int_0^1  \left(\sum_{i=1}^n t^{c_i-\frac{1}{m}}x_i
\right)^m dt.
\end{aligned}
$$
By the definition of Riemann integral, we have
$$
\mathcal{C} {\bf x}^m = \lim_{k \rightarrow \infty}\sum_{j=1}^k
\frac{\left(\sum_{i=1}^n (\frac{j}{k})^{c_i-\frac{1}{m}} x_i
\right)^m}{k}.
$$
Let $\mathcal{C}_k$ be the symmetric tensor such that
$$
\begin{aligned}
\mathcal{C}_k {\bf x}^m=&\sum_{j=1}^k \frac{\left(\sum_{i=1}^n (\frac{j}{k})^{c_i-\frac{1}{m}} x_i \right)^m}{k} \\
= & \sum_{j=1}^k \left(\sum_{i=1}^n \frac{(\frac{j}{k})^{c_i-\frac{1}{m}} }{k^{\frac{1}{m}}}x_i \right)^m \\
=& \sum_{j=1}^k \left(\langle u^j, {\bf x} \rangle \right)^m,
\end{aligned}
$$
where
$$
u^j=\left(\frac{(\frac{j}{k})^{c_1-\frac{1}{m}}
}{k^{\frac{1}{m}}},\cdots, \frac{(\frac{j}{k})^{c_n-\frac{1}{m}}
}{k^{\frac{1}{m}}}\right) \in \mathbb{R}^n, \ j=1,\cdots,k.
$$
Let ${\rm CD}_{m,n}$ denote the set consisting of all completely decomposable tensor with order $m$ and dimension $n$. From  \cite[Theorem 1]{LQX}, ${\rm CD}_{m,n}$ is a closed convex cone when $m$ is even. It then follows that
$\mathcal{C}=\lim_{k \rightarrow \infty}\mathcal{C}_k$ is also a  completely positive tensor.
\qed

\subsection{Even order symmetric weakly diagonally dominated tensors have SOS decompositions}
In this section, we establish that even order symmetric weakly diagonally dominated tensors have SOS decompositions.
Firstly, we give the definition of weakly diagonally dominated tensors. To do this, we introduce an index set $\Delta_{\mathcal{A}}$
associated with a tensor $\mathcal{A}$.
Now, let $\mathcal{A}$ be a tensor with order $m$ and dimension $n$, and let $f_{\mathcal{A}}$ be its associated homogeneous polynomial such that
$f_{\mathcal{A}}({\bf x})=\mathcal{A} {\bf x}^m$. We then define the index set $\Delta_{\mathcal{A}}$ as $\Delta_{f}$ with $f=f_{\mathcal{A}}$, as given as in
(\ref{e23}).

\begin{definition}\label{def31}
We say  $\mathcal{A}$ is a {\bf diagonally dominated tensor} if, for each $i=1,\cdots,n$,
\[
 a_{i i \cdots i} \ge \sum_{(i_2, \cdots, i_m) \neq (i \cdots i)} |a_{i i_2 \cdots i_m}|.
\]
We say $\mathcal{A}$ is a {\bf weakly diagonally dominated tensor} if, for each $i=1,\cdots,n$,
\[
a_{i i \cdots i} \ge \sum_{(i_2 \cdots i_m) \neq (i \cdots i), \atop {(i,i_2 \cdots,i_m) \in \Delta_{\mathcal{A}}}} |a_{i i_2 \cdots i_m}|.
\]
\end{definition}
Clearly, any diagonally dominated tensor is a weakly diagonally dominated tensor. However, the converse is, in general, not true.
\begin{theorem}\label{them31}
Let $\mathcal{A}$ be a symmetric weakly diagonally dominated tensor with order $m$ and dimension $n$. Suppose that $m$ is even. Then, $\mathcal{A}$ has an
SOS tensor decomposition.
\end{theorem}
\proof
Denote $I=\{(i,\cdots,i)\mid 1 \le i \le n\}$.
Let ${\bf x} \in \mathbb{R}^n$. Then,
\begin{eqnarray*}
 \mathcal{A}{\bf x}^m & = & \sum_{i=1}^n a_{i i \cdots i} x_i^m + \sum_{(i_1,\cdots,i_m) \notin I} a_{i_1 i_2 \cdots i_m} x_{i_1} x_{i_2} \cdots x_{i_m} \\
 & = & \sum_{i=1}^n \left(a_{i i \cdots i}-\sum_{(i_2 \cdots i_m) \neq (i \cdots i)\atop { (i,i_2 \cdots,i_m) \in \Delta_{\mathcal{A}}}} |a_{i i_2 \cdots i_m}| \right) x_i^m + \\
 & &  \sum_{i=1}^n \sum_{(i_2 \cdots i_m) \neq (i \cdots i) \atop { (i,i_2 \cdots,i_m) \in \Delta_{\mathcal{A}}}} |a_{i i_2 \cdots i_m}|  x_i^m + \sum_{(i_1,\cdots,i_m) \notin I} a_{i_1 i_2 \cdots i_m} x_{i_1} x_{i_2} \cdots x_{i_m} \\
 & = & \sum_{i=1}^n \left(a_{i i \cdots i}-\sum_{(i_2 \cdots i_m) \neq (i \cdots i) \atop (i,i_2 \cdots,i_m) \in \Delta_{\mathcal{A}}} |a_{i i_2 \cdots i_m}| \right) x_i^m  \\
 & &  +\sum_{i=1}^n \sum_{(i_2 \cdots i_m) \neq (i \cdots i) \atop { (i,i_2 \cdots,i_m) \in \Delta_{\mathcal{A}}}} |a_{i i_2 \cdots i_m}|  x_i^m +  \sum_{i=1}^n \sum_{(i_2 \cdots i_m) \neq (i \cdots i) \atop { (i,i_2 \cdots,i_m) \in \Delta_{\mathcal{A}}}}a_{i i_2 \cdots i_m} x_{i} x_{i_2} \cdots x_{i_m} \\
 & & + \sum_{i=1}^n \sum_{(i_2 \cdots i_m) \neq (i \cdots i) \atop { (i,i_2 \cdots,i_m) \notin \Delta_{\mathcal{A}}}}a_{i i_2 \cdots i_m} x_{i} x_{i_2} \cdots x_{i_m}
\end{eqnarray*}
Define
\[
h({\bf x})= \sum_{i=1}^n \sum_{(i_2 \cdots i_m) \neq (i \cdots i) \atop { (i,i_2 \cdots,i_m) \in \Delta_{\mathcal{A}}}} |a_{i i_2 \cdots i_m}|  x_i^m +  \sum_{i=1}^n \sum_{(i_2 \cdots i_m) \neq (i \cdots i) \atop { (i,i_2 \cdots,i_m) \in \Delta_{\mathcal{A}}}}a_{i i_2 \cdots i_m} x_{i} x_{i_2} \cdots x_{i_m}.
\]
We now show that $h$ is a sums-of-squares polynomial.

To see $h$ is indeed sums-of-squares, from Lemma \ref{lema21}, it suffices to show that
\[
\hat{h}({\bf x}):= \sum_{i=1}^n \sum_{(i_2 \cdots i_m) \neq (i \cdots i) \atop (i,i_2 \cdots,i_m) \in \Delta_{\mathcal{A}}} |a_{i i_2 \cdots i_m}|  x_i^m -  \sum_{i=1}^n \sum_{(i_2 \cdots i_m) \neq (i \cdots i) \atop (i,i_2 \cdots,i_m) \in \Delta_{\mathcal{A}}}|a_{i i_2 \cdots i_m}| x_{i} x_{i_2} \cdots x_{i_m}
\]
is a polynomial which always takes nonnegative values. As $\hat{h}$ is a homogeneous polynomial with degree $m$ on $\mathbb{R}^n$, let $\hat{\mathcal{H}}$ be a
symmetric tensor with order $m$ and dimension $n$ such that
$\hat{h}({\bf x})=\hat{\mathcal{H}}{\bf x}^m$.
Since $\mathcal{A}$ is symmetric, the nonzero entries of  $\hat{\mathcal{H}}$ are the same as the corresponding entries of $\mathcal{A}$.
Now, let $\lambda$ be an arbitrary $H$-eigenvalue of $\hat{\mathcal{H}}$, from the Gershgorin Theorem for eigenvalues of tensors \cite{Qi05}, we have
\[
\left|\lambda - \sum_{(i_2 \cdots i_m) \neq (i \cdots i) \atop (i,i_2 \cdots,i_m) \in \Delta_{\mathcal{A}}} |a_{i i_2 \cdots i_m}| \right| \le \sum_{(i_2 \cdots i_m) \neq (i \cdots i) \atop (i,i_2 \cdots,i_m) \in \Delta_{\mathcal{A}}} |a_{i i_2 \cdots i_m}|.
\]
So, we must have $\lambda \ge 0$. This shows that all $H$-eigenvalues of $\hat{\mathcal{H}}$ must be nonnegative, and so, $\hat{\mathcal{H}}$ is positive semi-definite \cite{Qi05}. Thus, $\hat{h}$ is
a polynomial which always takes nonnegative values.

Now, as $\mathcal{A}$ is a weakly diagonally dominated tensor and $m$ is even,
\[
\sum_{i=1}^n \left(a_{i i \cdots i}-\sum_{(i_2 \cdots i_m) \neq (i \cdots i) \atop (i,i_2 \cdots,i_m) \in \Delta_{\mathcal{A}}} |a_{i i_2 \cdots i_m}| \right) x_i^m
\]
is an SOS polynomial. Moreover, from the definition of $\Delta_{\mathcal{A}}$, for each $(i_1 \cdots i_m) \notin \Delta_{\mathcal{A}}$, $a_{i_1 \cdots i_m} \ge 0$ and
$x_{i_1} \cdots x_{i_m}$ is a squares term. Then,
\[
\sum_{i=1}^n \sum_{(i_2 \cdots i_m) \neq (i \cdots i) \atop { (i,i_2 \cdots,i_m) \notin \Delta_{\mathcal{A}}}}a_{i i_2 \cdots i_m} x_{i} x_{i_2} \cdots x_{i_m}
\]
is also a sums-of-square polynomial.
Thus, $\mathcal{A}$ has an SOS tensor decomposition. \qed

As a diagonally dominated tensor is weakly diagonally dominated, the following corollary follows immediately.
\begin{corollary}
Let $\mathcal{A}$ be a  symmetric diagonally dominated tensor with even order $m$ and dimension $n$. Then, $\mathcal{A}$ has an SOS tensor decomposition.
\end{corollary}

\subsection{The absolute tensor of an even order symmetric positive semi-definite $Z$-tensor has an SOS decomposition}
Let $\mathcal{A}$ be an order $m$ dimension $n$ tensor. If all off-diagonal elements of $\mathcal{A}$  are non-positive,
then $\mathcal{A}$ is called a $Z$-tensor \cite{Zhang12}. A $Z$-tensor $\mathcal{A}=(a_{i_1,\ldots,i_m})$  can be written as
\begin{equation}\label{e31} \mathcal{A}=\mathcal{D}-\mathcal{C}, \end{equation}
where $\mathcal{D}$ is a diagonal tensor where its $i$th diagonal elements equals $a_{i i \ldots i}$, $i=1,\ldots,n$, and $\mathcal{C}$ is a nonnegative tensor (or a tensor with nonnegative entries)  such that diagonal entries all equal to zero.
We now define the absolute tensor of $\mathcal{A}$ by
$$|\mathcal{A}|=|\mathcal{D}|+\mathcal{C}.$$


Note that all even order symmetric positive semi-definite $Z$-tensors have SOS decompositions \cite{Hu14,HLQS},
a natural interesting question would be: do all absolute tensors of even order symmetric positive semi-definite $Z$-tensors have SOS decompositions?
Below, we provide an answer for this question.

\bt\label{them34} Let $\mathcal{A}$ be a symmetric $Z$-tensor with even order $m$ and dimension $n$ defined as in (\ref{e31}).
If $\mathcal{A}$ is positive semi-definite, then $|\mathcal{A}|$ has an SOS tensor decomposition.
\et
\proof Let $\mathcal{A}=(a_{i_1 \ldots i_m})$ be a symmetric positive semi-definite $Z$-tensor. From (\ref{e31}), we have $\mathcal{A}= \mathcal{D}-\mathcal{C}$, where $\mathcal{D}$ is a diagonal tensor where the diagonal entries
of $\mathcal{D}$ is $d_i:=a_{i \ldots i}, i\in [n]$ and
$\mathcal{C}=(c_{i_1i_2\cdots i_m})$ is a nonnegative tensor with zero diagonal entries. Define three index sets as follows:
$$
\begin{aligned} I=&\{(i_1, i_2,\cdots,i_m)\in [n]^m~|~i_1=i_2=\cdots=i_m\};\\
\Omega=&\{(i_1, i_2,\cdots,i_m)\in [n]^m~|~ c_{i_1i_2\cdots i_m}\neq0~~ and ~~(i_1, i_2,\cdots,i_m) \notin I\}; \\
\Delta=&\{(i_1, i_2,\cdots,i_m)\in \Omega~|~ c_{i_1i_2\cdots i_m}>0~~{\bf or}~~{\rm at~~least~~one~~index~~in}~(i_1, i_2,\cdots,i_m)~~{\rm exists~~odd~~times}\}.
\end{aligned}
$$
Let $f({\bf x})=|\mathcal{A}|{\bf x}^m$ and define
  a polynomial $\hat{f}$ by
$$\hat{f}({\bf x})=\sum_{i=1}^nd_ix_i^m-\sum_{(i_1,i_2,\cdots,i_m)\in \Delta}|c_{i_1i_2\cdots i_m}|x_{i_1}x_{i_2}\cdots x_{i_m}.$$
 From Lemma \ref{lema21}, to see polynomial $f({\bf x})=|\mathcal{A}|{\bf x}^m$ is a sums-of-squares polynomial, we only need to show that
$\hat{f}$ always takes nonnegative value. To see this,
as $\mathcal{A}$ is positive semi-definite, we have $d_i\geq 0$.
Since $c_{i_1i_2\cdots i_m}\geq0$, $i_j\in [n],~j\in [m]$, it follows that
$$
\begin{aligned}
\hat{f}({\bf x})=&\sum_{i=1}^nd_ix_i^m-\sum_{(i_1,i_2,\cdots,i_m)\in \Delta}c_{i_1i_2\cdots i_m}x_{i_1}x_{i_2}\cdots x_{i_m}\\
=&\sum_{i=1}^nd_ix_i^m-\sum_{(i_1,i_2,\cdots,i_m)\in \Omega}c_{i_1i_2\cdots i_m}x_{i_1}x_{i_2}\cdots x_{i_m}+\sum_{(i_1,i_2,\cdots,i_m)\in \Omega \backslash \Delta}c_{i_1i_2\cdots i_m}x_{i_1}x_{i_2}\cdots x_{i_m}\\
=&\mathcal{A}{\bf x}^m+\sum_{(i_1,i_2,\cdots,i_m)\in \Omega \backslash \Delta}c_{i_1i_2\cdots i_m}x_{i_1}x_{i_2}\cdots x_{i_m} \\
\geq&0.
\end{aligned}
$$
Here, the last inequality follows from the fact that $m$ is even, $\mathcal{A}$ is positive semi-definite and
$x_{i_1}x_{i_2}\cdots x_{i_m}$ is a square term if $(i_1,i_2,\cdots,i_m)\in \Omega \backslash \Delta$.
Thus, the desired result follows.
\qed

\subsection{SOS tensor decomposition for even order symmetric extended $Z$-tensors}
In this subsection, we introduce a new class of symmetric tensor which extends symmetric $Z$-tensors to the cases where the off-diagonal elements can be positive, and examine its SOS tensor decomposition.

Let $f$ be a polynomial on $\mathbb{R}^n$ with degree $m$.  Let $f_{m,i}$ be the coefficient of $f$ associated with $x_i^{m}$, $i\in [n]$. We say $f$ is an extended $Z$-polynomial if there exist $s \in \mathbb{N}$ with $s \le n$ and index sets $\Gamma_l \subseteq \{1,\cdots,n\}$, $l=1,\cdots,s$ with  $\bigcup_{l=1}^s \Gamma_l=\{1,\cdots,n\}$ and $\Gamma_{l_1} \cap \Gamma_{l_2} =\emptyset$ for all $l_1 \neq l_2$  such that
\[
f({\bf x})=\sum_{i=1}^n f_{m,i} x_i^{m}+\sum_{l=1}^s \sum_{\alpha_l \in \Omega_l }f_{\alpha_l} {\bf x}^{\alpha_l},\]
where
\[
\Omega_l=\left\{\alpha\in ([n]\cup \{0\})^n:~|\alpha|=m, {\bf x}^{\alpha}=x_{i_1}x_{i_2}\cdots x_{i_m}, \{i_1,\cdots,i_m\}\subseteq \Gamma_l,~
\mbox{and} ~~ \alpha \neq m {\bf e_i}, \ i=1,\cdots,n \right \}
\]
for each $l=1,\cdots,s$ and either one of the following two conditions holds:
\begin{itemize}
\item[{\rm (1)}] $f_{\alpha_l}=0$ for all but one $\alpha_l \in \Omega_l$;
\item[{\rm (2)}] $f_{\alpha_l} \le 0$ for all $\alpha_l \in \Omega_l$.
\end{itemize}
We now say a symmetric tensor $\mathcal{A}$ is an {\bf extended $Z$-tensor} if its associated polynomial $f_{\mathcal{A}}({\bf x})=\mathcal{A}{\bf x}^m$ is a
an extended $Z$-polynomial.

From the definition, it is clear that any $Z$-tensor is an extended $Z$-tensor with $s=1$ and $\Gamma_1=\{1,\cdots,n\}$. On the other hand, an extended $Z$-tensor allows a few elements of the off-diagonal elements to be positive, and so, an extended $Z$-tensor need not to be a $Z$-tensor. For example, consider a symmetric tensor $\mathcal{A}$ where its associated polynomial   $f_{\mathcal{A}}({\bf x})=\mathcal{A}{\bf x}^m=x_1^6+x_2^6+x_3^6+x_4^6+4x_1^3x_2^3+6x_3^2x_4^4$. It can be easily see that $\mathcal{A}$ is an extended $Z$-tensor but not a $Z$-tensor (as there are positive off-diagonal elements).    In \cite{LQX15}, partially $Z$-tensors are introduced.   There is no direct relation between these two concepts, except that both of them contain $Z$-tensors. But they do have intersection which is larger than the set of all $Z$-tensors.   Actually, the example just discussed is not a $Z$-tensor, but it is an extended $Z$-tensor and a partially $Z$-tensor as well.

We now see that any positive semi-definite extended $Z$-tensor has an SOS tensor decomposition. To achieve this, we recall the following useful lemma, which provides us a
simple criterion for determining whether a homogeneous polynomial with only one mixed term is a sum of
squares polynomial or not.

\bl\label{lema52}$^{\cite{FK11}}$ Let $b_1,b_2,\cdots,b_n\geq 0$ and $d\in \mathbb{N}$.
Let $a_1,a_2,\cdots,a_n\in \mathbb{N}$ be such that $\sum_{i=1}^n=2d$. Consider the homogeneous
polynomial $f({\bf x})$ defined by
$$f({\bf x})=b_1x_1^{2d}+\cdots+b_nx_n^{2d}-\mu x_1^{a_1}\cdots x_n^{a_n}.$$
Let $\mu_0=2d\prod _{a_i\neq 0, 1\leq i\leq n}(\frac{b_i}{a_i})^{\frac{a_i}{2d}}.$
Then, the following statements are equivalent:
\begin{itemize}
 \item[{\rm (i)}] $f$ is a nonnegative polynomial i.e. $f({\bf x})\geq 0$ for all ${\bf x}\in \mathbb{R}^n$;
 \item[{\rm (ii)}] either $|\mu|\leq \mu_0$ or $\mu<\mu_0$ and all $a_i$ are even;
 \item[{\rm (iii)}] $f$ is an SOS polynomial.
\end{itemize}

\bt \label{th:5.2}
Let $\mathcal{A}$ be an even order positive semi-definite extended $Z$-tensor. Then, $\mathcal{A}$ has an SOS tensor decomposition.
\et
\proof
Let $f_{\mathcal{A}}({\bf x})=\mathcal{A}{\bf x}^m$. As $\mathcal{A}$ is a positive semi-definite symmetric extended $Z$-tensor, there exist $s \in \mathbb{N}$ and index sets $\Gamma_l \subseteq \{1,\cdots,n\}$, $l=1,\cdots,s$ with  $\bigcup_{l=1}^s \Gamma_l=\{1,\cdots,n\}$ and $\Gamma_{l_1} \cap \Gamma_{l_2} =\emptyset$ for all $l_1 \neq l_2$  such that
for all ${\bf x} \in \mathbb{R}^n$
\[
f({\bf x})=\sum_{i=1}^n f_{m,i} x_i^{m}+\sum_{l=1}^s \sum_{\alpha_l \in \Omega_l}f_{\alpha_l} {\bf x}^{\alpha_l}
\]
such that, for each $l=1,\cdots,s$, either one of the following two condition holds:
(1) $f_{\alpha_l}=0$ for all but one $\alpha_l \in \Omega_l$; (2) $f_{\alpha_l} \le 0$ for all $\alpha_l \in \Omega_l$.
Define, for each $l=1,\cdots,s$,
\[
 h_l({\bf x}):=\sum_{i \in \Gamma_l} f_{m,i} x_i^{m}+ \sum_{\alpha_l \in \Omega_l}f_{\alpha_l} {\bf x}^{\alpha_l}.
\]
It follows that each $h_l$ is an extended $Z$-polynomial. Moreover, from the construction, $\sum_{l=1}^s h_l=f_{\mathcal{A}}$ and so, $\inf_{{\bf x} \in \mathbb{R}^n} \sum_{l=1}^s h_l({\bf x})=0$. Note that each $h_l$ is also a homogeneous polynomial, and hence
$\inf_{{\bf x} \in \mathbb{R}^n} h_l({\bf x}) \le 0$. Noting that each $h_l$ is indeed a polynomial on $(x_i)_{i \in \Gamma_l}$, $\bigcup_{l=1}^s \Gamma_l=\{1,\cdots,n\}$ and $\Gamma_{l_1} \cap \Gamma_{l_2} =\emptyset$ for all $l_1 \neq l_2$, we have  $\inf_{{\bf x} \in \mathbb{R}^n} \sum_{l=1}^s h_l({\bf x}) = \sum_{l=1}^s \inf_{{\bf x} \in \mathbb{R}^n} h_l({\bf x})$. This enforces that $\inf_{{\bf x} \in \mathbb{R}^n} h_l({\bf x})=0$. In particular,
each $h_l$ is a polynomial which takes nonnegative values. We now see that $h_l$, $1 \le l \le s$, are SOS polynomial. Indeed, if $f_{\alpha_l}=0$ for all but one $\alpha_l \in \Omega_l$, then $h_l$ is a homogeneous polynomial with only a mixed term, and so,
 Lemma \ref{lema52} implies that $h_l$ is a SOS polynomial. On the other hand, if $f_{\alpha_l} \le 0$ for all $\alpha_l \in \Omega_l$, $h_l$ corresponds to a $Z$-tensor, and so, $h_l$ is also a SOS polynomial in this case because any positive semi-definite $Z$-tensor has an SOS tensor decomposition \cite{Hu14}.
Thus, $f_{\mathcal{A}}=\sum_{l=1}^sh_l$ is also a SOS polynomial, and hence the conclusion follows. \qed

\begin{Remark}\label{remark:3.3}
A close inspection of the above proof indicates that we indeed shows that the associated polynomial $f_{\mathcal{A}}({\bf x})=\mathcal{A}{\bf x}^m$ satisfies
$f_{\mathcal{A}}=\sum_{l=1}^s h_l$ where each $h_l$ is an SOS polynomial in $(x_i)_{i \in \Gamma_l}$.
\end{Remark}

\subsection{Even order symmetric $B_0$-tensors have SOS decompositions}
In this part, we show that even order symmetric $B_0$ tensors have SOS tensor decompositions.
Recall that a tensor  $\mathcal{A}=(a_{i_1i_2\cdots i_m})$ with order $m$ and dimension $n$ is called a $B_0$-tensor \cite{qi14} if
 \[
  \sum_{i_2, \cdots, i_m=1}^n a_{i i_2 \cdots i_m} \ge 0
 \]
and
\[
\frac{1}{n^{m-1}} \sum_{i_2, \cdots, i_m=1}^n a_{i i_2 \cdots i_m} \ge a_{i j_2 \cdots j_m} \mbox{ for all } (j_2,\cdots, j_m) \neq (i, \cdots, i).
\]
To establish that a $B_0$-tensor has an SOS tensor decomposition, we first present the SOS tensor decomposition of the all-one-tensor. We say $\mathcal{E}$ is an all-one-tensor if
 with each of its elements of $\mathcal{E}$ is equal to one.
\bl\label{lema32} Let $\mathcal{E}$ be  an even order all-one-tensor. Then, $\mathcal{E}$ has an SOS tensor decomposition.
\el
\proof Let $\mathcal{E}=(e_{i_1i_2\cdots i_m})$ be an all-one-tensor with even order $m$ and dimension $n$.
For all ${\bf x}\in \mathbb{R}^n$, one has
$$
\begin{aligned}
\mathcal{E}{\bf x}^m=&\sum_{i_1,i_2,\cdots,i_m \in [n]}e_{i_1i_2\cdots i_m}x_{i_1}x_{i_2}\cdots x_{i_m}\\
=&\sum_{i_1,i_2,\cdots,i_m \in [n]}x_{i_1}x_{i_2}\cdots x_{i_m} \\
=&(x_1+x_2+\cdots+x_n)^m \\
\geq & 0,
\end{aligned}
$$
which implies that $\mathcal{E}$ has an SOS tensor decomposition. \qed
Let $J\subset [n]$. $\mathcal{E}^J$ is called a partially all-one-tensor if its elements are defined such that
$e_{i_1i_2\cdots i_m}=1$, $i_1, i_2, \cdots, i_m\in J$ and $e_{i_1i_2\cdots i_m}=0$ for the others.
{Similar to Lemma 3.2}, it is easy to check that all even order partially all-one-tensors
have SOS decompositions.

We also need the following characterization of $B_0$-tensors established in \cite{qi14}.
\begin{lemma}\label{lema33}
Suppose that $\mathcal{A}$ is a $B_0$-tensor with order $m$ and dimension $n$.
Then either $\mathcal{A}$ is a diagonally dominated symmetric $M$-tensor itself, or we have
 \[
 \mathcal{A}=\mathcal{M}+ \sum_{k=1}^s h_k \mathcal{E}^{J_k},
 \]
where $\mathcal{M}$ is a diagonally dominated symmetric $M$-tensor, $s$ is a positive integer, $h_k>0$ and $J_k  \subseteq \{1,\cdots,n\}$, for $k=1,\cdots,s$, and $J_k \cap J_l=\emptyset$,
for $k \neq l$.
\end{lemma}

From Theorem \ref{them31}, Lemma \ref{lema32} and Lemma \ref{lema33}, we have the following result.

\bt\label{them32} All even order symmetric $B_0$-tensors have SOS tensor decompositions.
\et

{
{ Before we move on to the next part,} we note that, stimulated by $B_0$-tensors in \cite{qi14},
symmetric double $B$-tensors, symmetric quasi-double $B_0$-tensors and symmetric $MB_0$-tensors {have been} studied
in \cite{LCL,LCQL}. Below, we briefly explain that, using a similar method of proof as above, these three classes of tensors all have SOS decompositions. To do this, let us
recall the definitions of these three classes of tensors. }

For a real symmetric tensor $\mathcal{B}=(b_{i_1i_2\cdots i_m})$ with order $m$ and dimension $n$,
denote $$\beta_i(\mathcal{B})=\max_{j_2,\cdots,j_m\in [n], (i,j_2,\cdots,j_m)\notin I}\{0, b_{ij_2\cdots j_m}\};$$
$$ \Delta_i(\mathcal{B})=\sum_{j_2,\cdots,j_m\in [n],(i,j_2,\cdots,j_m)\notin I}(\beta_i(\mathcal{B})-b_{ij_2\cdots j_m});$$
$$\Delta^i_j(\mathcal{B})=\Delta_j(\mathcal{B})-(\beta_j(\mathcal{B})-b_{jii\cdots i}),~i\neq j.$$
As defined in  \cite[Definition 3]{LCL},  $\mathcal{B}$ is called a double $B$-tensor if, $b_{ii\cdots i}> \beta_i(\mathcal{B}),$ for all $i\in [n]$ and for all
$i, j\in [n], i\neq j$ such that
$$b_{ii\cdots i}-\beta_i(\mathcal{B})\geq \Delta_i(\mathcal{B})$$
and
$$(b_{ii\cdots i}-\beta_i(\mathcal{B}))(b_{jj\cdots j}-\beta_j(\mathcal{B}))>\Delta_i(\mathcal{B})\Delta_j(\mathcal{B}).$$
If $b_{ii\cdots i}> \beta_i(\mathcal{B}),$ for all $i\in [n]$ and
$$(b_{ii\cdots i}-\beta_i(\mathcal{B}))(b_{jj\cdots j}-\beta_j(\mathcal{B})-\Delta^i_j(\mathcal{B}))\geq (\beta_j(\mathcal{B})-b_{ji\cdots i})\Delta_i(\mathcal{B}),$$
then tensor $\mathcal{B}$ is called a quasi-double $B_0$-tensor (see Definition 2 of \cite{LCQL}).

Let $\mathcal{A}=(a_{i_1i_2\cdots i_m})$ such that
$$a_{i_1i_2\cdots i_m}=b_{i_1i_2\cdots i_m}-\beta_{i_1}(\mathcal{B}),~{\rm for~~ all}~~ i_1\in [n].$$
If $\mathcal{A}$ is an $M$-tensor, then $B$ is called an $MB_0$-tensor (see Definition 3 of \cite{LCQL}).
It was shown in \cite{LCQL} that all quasi-double $B_0$-tensors are $MB_0$-tensors.

 In \cite{LCL}, Li et al. proved that, for any symmetric double $B$-tensor $\mathcal{B}$, either
$\mathcal{B}$ is a doubly strictly diagonally dominated (DSDD) $Z$-tensor, or $\mathcal{B}$ can be decomposed to
the sum of a DSDD $Z$-tensor and several positive multiples of partially all-one-tensors (see Theorem 6 of \cite{LCL}).
From Theorem 4 of \cite{LCL}, we know that an even order symmetric DSDD $Z$-tensor is positive definite. This together with
the fact that any positive semi-definite $Z$-tensor has an SOS tensor decomposition \cite{Hu14} implies that any even order symmetric double $B$-tensor $\mathcal{B}$
has an SOS tensor decomposition. Moreover, from Theorem 7 of \cite{LCQL}, we know that, for any symmetric $MB_0$-tensor, it is either an $M$-tensor itself or it can be decomposed
as the sum of an $M$-tensor and several positive multiples of partially {all-one-tensors}. As even order symmetric $M$-tensors are positive semi-definite $Z$-tensors \cite{Zhang12} which have, in particular, SOS decomposition, we see that any even order symmetric $MB_0$ tensor also has an SOS tensor decomposition. Combining these and noting that any quasi-double $B_0$-tensor is an $MB_0$-tensor, we arrive at the following conclusion.

\bt\label{them33} Even order symmetric double $B$-tensors, even order symmetric quasi double $B_0$-tensors
and even order symmetric $MB_0$-tensors all have SOS tensor decompositions.
\et

{ \subsection{Even order symmetric $H$-tensors with nonnegative diagonal elements have SOS decompositions}

In this part, we show that any even order symmetric $H$-tensor with nonnegative diagonal elements has an SOS tensor decomposition.
{Recall that  an $m$th order $n$ dimensional tensor $\mathcal{A}=(a_{i_1i_2\cdots i_m})$, it's comparison tensor is
defined by $\mathcal{M}(\mathcal{A})=(m_{i_1i_2\cdots i_m})$ such that
$$m_{ii\cdots i}=|a_{ii\cdots i}|,~~{\rm and}~~ m_{i_1i_2\cdots i_m}=-|a_{i_1i_2\cdots i_m}|,
$$
for all $i,i_1,\cdots,i_m\in [n], (i_1,i_2,\cdots,i_m)\notin I$.
Then, tensor $\mathcal{A}$ is called an $H$-tensor \cite{Ding13} if there exists a tensor $\mathcal{Z}$ with nonnegative entries such that
$\mathcal{M}(\mathcal{A})=s \mathcal{I}- \mathcal{Z}$ and $s \ge \rho(\mathcal{Z})$, where $\mathcal{I}$ is the identity tensor and $\rho(\mathcal{Z})$ is the spectral
radius of $\mathcal{Z}$ defined as the maximum of modulus of all eigenvalues of $\mathcal{Z}$. If $s > \rho(\mathcal{Z})$, then $\mathcal{A}$ is called
a nonsingular $H$-tensor.} A characterization for nonsingular $H$-tensors was given in \cite{Ding13} which states
$\mathcal{A}$ is a nonsingular
$H$-tensor if and only if there exists an enteritis positive vector ${\bf y}=(y_1,y_2,\cdots,y_n)\in \mathbb{R}^n$ such that
$$|a_{ii\cdots i}|y_i^{m-1}>\sum_{(i,i_2,\cdots,i_m)\notin I}|a_{ii_2\cdots i_m}|y_{i_2}y_{i_3}\cdots y_{i_m},~~\forall~i\in[n].$$
We note that the above definitions were first introduced in \cite{Ding13}. These were further
 examined in  \cite{Kannan15,LWZ} where the authors in \cite{LWZ} { referred} nonsingular $H$-tensors simply as $H$-tensors and the authors in \cite{Kannan15} { referred} nonsingular $H$-tensors as strong $H$-tensors.

\bt\label{them35} Let $\mathcal{A}=(a_{i_1i_2\cdots i_m})$ be a symmetric $H$-tensor with even order $m$ dimension $n$. Suppose that all the diagonal elements of $\mathcal{A}$ are nonnegative. Then, $\mathcal{A}$ has an SOS tensor decomposition.
\et
\proof We first show that any nonsingular $H$-tensor with positive diagonal elements has an SOS tensor decomposition. Let $\mathcal{A}=(a_{i_1i_2\cdots i_m})$ be a nonsingular  $H$-tensor with even order $m$ dimension $n$ such that $a_{i i \cdots i} > 0$, $i \in [n]$.
Then, there exists a vector ${\bf y}=(y_1,\cdots,y_n)^T \in \mathbb{R}^n$ with $y_i>0$, $i=1,\cdots,n$, such that
\begin{equation}\label{e34}
a_{ii\cdots i}y_i^{m-1} > \sum_{(i,i_2,\cdots,i_m)\notin I}|a_{ii_2\cdots i_m}|y_{i_2}y_{i_3}\cdots y_{i_m},~\forall~i\in [n].
\end{equation}
To prove the conclusion, by Lemma \ref{lema21}, we only need to prove
$$\hat{f}_{\mathcal{A}}({\bf x})=\sum_{i\in [n]}a_{ii\cdots i}x_i^m-\sum_{(i_1,i_2,\cdots,i_m)\in \Delta_{\mathcal{A}}}|a_{i_1i_2\cdots i_m}|x_{i_1}x_{i_2}\cdots x_{i_m}\geq 0,~\forall~{\bf x}\in \mathbb{R}^n.$$
From (\ref{e34}), we know that
\begin{equation}\label{e35}
\begin{aligned}
\hat{f}_{\mathcal{A}}({\bf x})\geq &\sum_{i\in [n]} \left(\sum_{(i,i_2,\cdots,i_m)\notin I}|a_{ii_2\cdots i_m}|y_i^{1-m}y_{i_2}y_{i_3}\cdots y_{i_m} x_i^m\right)
-\sum_{(i_1,i_2,\cdots,i_m)\in \Delta_{\mathcal{A}}}|a_{i_1i_2\cdots i_m}|x_{i_1}x_{i_2}\cdots x_{i_m}.
\end{aligned}
\end{equation}

 Here, for any fixed tuple $(i_1^0,i_2^0,\cdots,i_m^0)\in \Delta_{\mathcal{A}}$,
assume $(i_1^0,i_2^0,\cdots,i_m^0)$ is constituted by $k$ distinct indices $j^0_1,j^0_2,\cdots, j^0_k$, $k\leq m$,
which appear $s_1,s_2,\cdots,s_k$ times in $(i_1^0,i_2^0,\cdots,i_m^0)$ respectively, $s_l\in[m], l\in [k]$. Then, one has $s_1+s_2+\cdots+s_k=m$.
 Without loss of generality, we denote $a=|a_{i^0_1i^0_2\cdots i^0_m}|>0$. {  Let $\pi(i_1^0,i_2^0,\cdots,i_m^0)$ be the set consisting of all permutations of $(i_1^0,i_2^0,\cdots,i_m^0).$}
So, on the right side of (\ref{e35}), there are some terms corresponding to the fixed tuple $(i^0_1,i^0_2,\cdots,i^0_m)$ such that{
\begin{eqnarray*}
& & \sum_{(j^0_1,i_2,\cdots,i_m)\in \pi(i_1^0,i_2^0,\cdots,i_m^0)}|a_{j^0_1i_2\cdots i_m}|y_{j_1^0}^{1-m}y_{i_2}y_{i_3}\cdots y_{i_m} x_{j_1^0}^m\\
& & +\sum_{(j^0_2,i_2,\cdots,i_m)\in \pi(i_1^0,i_2^0,\cdots,i_m^0)}|a_{j^0_2i_2\cdots i_m}|y_{j_2^0}^{1-m}y_{i_2}y_{i_3}\cdots y_{i_m} x_{j_2^0}^m\\
& & +\cdots\\
& & +\sum_{(j^0_k,i_2,\cdots,i_m)\in \pi(i_1^0,i_2^0,\cdots,i_m^0)}|a_{j^0_ki_2\cdots i_m}|y_{j_k^0}^{1-m}y_{i_2}y_{i_3}\cdots y_{i_m} x_{j_k^0}^m\\
& & -\sum_{(i_1,i_2,\cdots,i_m)\in \pi(i_1^0,i_2^0,\cdots,i_m^0)}|a_{i_1i_2\cdots i_m}|x_{i_1}x_{i_2}\cdots x_{i_m}\\
&= & \binom {m-1} {s_1-1} \binom {m-s_1} {s_2} \binom {m-s_1-s_2}{s_3}\cdots \binom {m-s_1-s_2\cdots -s_{k-1}} {s_k}ay_{j^0_1}^{s_1-m}y_{j^0_2}^{s_2}\cdots y_{j^0_k}^{s_k}x_{j^0_1}^m \\
& & +\binom {m-1} {s_2-1} \binom {m-s_2} {s_1} \binom {m-s_1-s_2}{s_3}\cdots \binom {m-s_1-s_2\cdots -s_{k-1}} {s_k}ay_{j^0_2}^{s_2-m}y_{j^0_1}^{s_1}y_{j^0_3}^{s_3}\cdots y_{j^0_k}^{s_k}x_{j^0_2}^m \\
& & +\cdots \cdots \\
& & +\binom {m-1} {s_k-1} \binom {m-s_k} {s_1} \binom {m-s_k-s_1}{s_2}\cdots \binom {m-s_k-s_1\cdots -s_{k-2}} {s_{k-1}}ay_{j^0_k}^{s_k-m}y_{j^0_1}^{s_1}\cdots y_{j^0_{k-1}}^{s_{k-1}}x_{j^0_k}^m \\
& & -\binom {m} {s_1} \binom {m-s_1} {s_2} \binom {m-s_1-s_2}{s_3}\cdots \binom {m-s_1-s_2\cdots -s_{k-1}} {s_k}ax_{j^0_1}^{s_1}x_{j^0_2}^{s_2} \cdots x_{j^0_k}^{s_k}\\
&  = & \frac{(m-1)!ay_{j^0_1}^{s_1}y_{j^0_2}^{s_2}\cdots y_{j^0_k}^{s_k}}{s_1!s_2!\cdots s_k!} \left[s_1\left(\frac{x_{j^0_1}}{y_{j^0_1}}\right)^m +s_2\left(\frac{x_{j^0_2}}{y_{j^0_2}}\right)^m +\cdots+s_k\left(\frac{x_{j^0_k}}{y_{j^0_k}}\right)^m -m\left(\frac{x_{j^0_1}}{y_{j^0_1}}\right)^{s_1} \left(\frac{x_{j^0_2}}{y_{j^0_2}}\right)^{s_2}\cdots \left(\frac{x_{j^0_k}}{y_{j^0_k}}\right)^{s_k}\right] \\
& \geq &  0,
\end{eqnarray*}}
where the last inequality follows the arithmetic-geometric inequality and the fact ${\bf y}>{\bf 0}$.
Thus, each tuple $(i_1,i_2,\cdots,i_m)\in \Delta_{\mathcal{A}}$ corresponds to a nonnegative value on the right side of (\ref{e35}), which implies that
$\hat{f}({\bf x})\geq 0$ for all ${\bf x}\in \mathbb{R}^n$. Hence, by Lemma \ref{lema21}, $\mathcal{A}$ has an SOS tensor decomposition.

Now, let $\mathcal{A}$ be a general $H$-tensor with nonnegative diagonal elements. Then, for each $\epsilon>0$, $\mathcal{A}_{\epsilon}:=\mathcal{A}+\epsilon \mathcal{I}$ is a
nonsingular $H$-tensor with positive diagonal elements. Thus, $\mathcal{A}_{\epsilon} \rightarrow \mathcal{A}$, and for each $\epsilon>0$, $\mathcal{A}_{\epsilon}$ has an SOS tensor decomposition. As ${\rm SOS}_{m,n}$ is a closed convex cone, we see that $\mathcal{A}$ also has an SOS tensor
decomposition and the desired results follows. \qed}

In \cite{DQW1, LQX}, SOS decomposition of some classes of Hankel tensors was given.

\setcounter{equation}{0}
\section{ The SOS-Rank of SOS tensor Decomposition}
In this section, we study the SOS-rank of SOS tensor decomposition. Let us formally define the SOS-rank of SOS tensor decomposition as follows.
Let $\mathcal{A}$ be a tensor with even order $m$ and dimension $n$. Suppose $\mathcal{A}$ has a SOS tensor decomposition.
As shown in Proposition 3.1, ${\rm SOS}_{m,n}={\rm SOS}_{m,n}^h$ where ${\rm SOS}_{m,n}$ is the SOS tensor cone and
${\rm SOS}_{m,n}^h$ is the cone consisting of all $m$th-order $n$-dimensional symmetric tensors such that
 $f_{\mathcal{A}}({\bf x}):=\langle \mathcal{A}, {\bf x}^ m \rangle$ is a polynomial which can be written as sums of finitely many \emph{homogeneous polynomials}.   Thus, there exists
$r\in \mathbb{N}$ such that the homogeneous polynomial
$f_{\mathcal{A}}({\bf x})=\mathcal{A}{\bf x}^m$ can be decomposed by
$$ f_{\mathcal{A}}({\bf x})=f_1^2({\bf x})+f_2^2({\bf x})+\cdots+f_r^2({\bf x}),~\forall ~{\bf x}\in \mathbb{R}^n,$$
where $f_i({\bf x}),~i\in [r]$ are homogeneous polynomials with degree $\frac{m}{2}$.
The minimum value $r$ is called the {\bf SOS-rank} of $\mathcal{A}$, and is denoted by
${\rm SOSrank}(\mathcal{A}).$

Let $C$ be a convex cone in the SOS tensor cone, that is, $C \subseteq {\rm SOS}_{m,n}$. We define the {\bf SOS-width} of the convex cone $C$  by
\[
\mbox {SOS-width}(C)=\sup\{ {\rm SOSrank}(\mathcal{A}): \mathcal{A} \in C\}.
\]
Here, we do not care about the minimum of the SOS-rank of all the possible tensors in the cone $C$ as it will be always zero.
 Recall that it was shown by Choi et al. in \cite[Theorem 4.4]{Ch95} that, an SOS homogeneous polynomial can be decomposed as sums of at most $\Lambda$
 many squares of homogeneous polynomials where \begin{equation}\label{eq:lambda}\Lambda=\frac{\sqrt{1+8a}-1}{2} \mbox{ and } a=\binom {n+m-1} {m}.\end{equation} This immediately gives us that

\bp \label{proposition} Let $\mathcal{A}$ be a tensor with even order $m$ and dimension $n$, $m, n \in \mathbb{N}$.
Suppose $\mathcal{A}$ has an SOS tensor decomposition.
Then, its SOS-rank satisfies
${\rm SOSrank(\mathcal{A})}\leq \Lambda$,
where $\Lambda$ is given in (\ref{eq:lambda}).
In particular, {\rm SOS-width}$({\rm SOS}_{m,n}) \leq \Lambda.$
\ep

In the matrix case, that is, $m=2$, the upper bound $\Lambda$ equals the dimension $n$ of the symmetric tensor which is tight in this case.
On the other hand, in general, the upper bound is of the order $n^{m/2}$ and need not to be tight.
However, for a class of structured tensors with bounded exponent (BD-tensors) that have SOS decompositions, we show that their SOS-rank is less or
equal to the dimension $n$ which is significantly smaller than the
upper bound in the above proposition. Moreover, in this case,
the SOS-width of the associated BD-tensor cone can be determined explicitly. To do this, let us recall the definition of polynomials with bounded exponent and define
the BD-tensors. Let $e \in \mathbb{N}$. Recall that $f$ is said to be a degree $m$ homogeneous polynomials on $\mathbb{R}^n$ with bounded exponent $e$ if

\[
f({\bf x})=\sum_{\alpha} f_{\alpha} {\bf x}^{\alpha}= \sum_{\alpha} f_{\alpha} x_1^{\alpha_1}\cdots x_n^{\alpha_n},
\]
where $0 \le \alpha_j \le e$ and $\sum_{j=1}^n \alpha_j=m$. We note that  degree $4$ homogeneous polynomials on $\mathbb{R}^n$ with bounded exponent $2$ is nothing but the bi-quadratic
forms in dimension $n$. Let us denote ${\rm BD}_{m,n}^{e}$ to be the set consists of all degree $m$ homogeneous polynomials on $\mathbb{R}^n$ with bounded exponent $e$.

An interesting result for characterizing when a positive semi-definite (PDF) homogeneous polynomial with bounded exponent has SOS tensor decomposition was established in \cite{Ch77} and can be stated as follows.
\begin{lemma}
Let $n \in \mathbb{N}$ with $n \ge 3$. Suppose $e,m$ are even numbers and $m \ge 4$.
\begin{itemize}
                                               \item[{\rm (1)}] If $n \ge 4$, then ${\rm BD}_{m,n}^{e} \cap {\rm PSD}_{m,n} \subseteq {\rm SOS}_{m,n}$ if and only if $m \ge en-2$;
                                               \item[{\rm (2)}] If $n =3$, then ${\rm BD}_{m,n}^{e} \cap {\rm PSD}_{m,n} \subseteq {\rm SOS}_{m,n}$ if and only if $m=4$ or $m \ge 3e-4$.
                                              \end{itemize}

\end{lemma}

Now, we say a symmetric tensor $\mathcal{A}$ is a BD-tensor with order $m$, dimension $n$ and exponent $e$ if $f({\bf x})=\mathcal{A}{\bf x}^m$ is a degree $m$ homogeneous polynomial on $\mathbb{R}^n$ with bounded exponent $e$.
We also define  ${\rm BD}_{m,n}$ to be the set consisting of all symmetric BD-tensors with order $m$, dimension $n$ and exponent $e$. It is clear that ${\rm BD}_{m,n}^e$ is a convex cone.

\bt
Let $n \in \mathbb{N}$ with $n \ge 3$. Suppose $e,m$ are even numbers and $m \ge 4$. Let $\mathcal{A}$ be a BD-tensor with order $m$, dimension $n$ and exponent $e$. Suppose that $\mathcal{A}$
has an SOS tensor decomposition. Then, we have ${\rm SOSrank} (\mathcal{A}) \le n$. Moreover, we have
$$\mbox{\rm SOS-width}({\rm BD}_{m,n}^{e} \cap {\rm SOS}_{m,n})=\left\{\begin{array}{ll}
1 & \mbox{ if } m=en \\
n & \mbox{ otherwise}.
\end{array}
\right.$$
\et

\begin{proof} As $\mathcal{A}$ is a BD-tensor and it has SOS decomposition, the preceding lemma implies that either {\rm (i)} $n \ge 4$ and $m \ge en-2$ {\rm (ii)} $n=3$ and $m=4$
and {\rm (iii)} $n=3$ and $m \ge 3e-4$. We now divide the discussion into these three cases.

Suppose that Case {\rm (i)} holds, i.e., $n \ge 4$ and $m \ge en-2$. From the construction, we have $m \le en$. If $m=en$, then $\mathcal{A}$ has the form $a x_1^{e}\cdots x_n^{e}$. Here,  $a \ge 0$ because
$\mathcal{A}$ has SOS decomposition and $e$ is an even number. In this case, ${\rm SOSrank}(\mathcal{A})=1$. Now, let $m=en-2$. Then,
\[
\mathcal{A} {\bf x}^m= x_1^{e} \cdots x_n^{e} \left(\sum_{(i,j) \in F} a_{ij} x_i^{-1} x_j^{-1}\right),
\]
for some $a_{ij} \in \mathbb{R}$, $(i,j)\in F$ and for some $F \subseteq \{1,\cdots,n\} \times \{1,\cdots,n\}$. As $e$ is an even number and $\mathcal{A}$ has SOS decomposition, we have
\[
\sum_{(i,j) \in F} a_{ij} x_i^{-1} x_j^{-1} \ge 0 \mbox{ for all } x_i \neq 0 \mbox{ and } x_j \neq 0.
\]
Thus, by continuity, $Q(t_1,\cdots,t_n)=\sum_{(i,j) \in F} a_{ij} t_i t_j$ is a positive semi-definite quadratic form, and so, is at most  sums of  $n$ many squares of linear functions in $t_1,\cdots,t_n$.
Let $Q(t_1,\cdots,t_n)=\sum_{k=1}^n \big[q_k(t_1,\cdots,t_n)\big]^2$ where $q_k$ are linear functions. Then,
\[
\mathcal{A} {\bf x}^m= x_1^{e} \cdots x_n^{e} \left(\sum_{i=1}^n \left[q_k(x_1^{-1},\cdots,x_n^{-1})\right]^2\right)=\sum_{i=1}^n\left(x_1^{e} \cdots x_n^{e}\left[q_k(x_1^{-1},\cdots,x_n^{-1})\right]^2\right),
\]
Note that $$x_1^{e} \cdots x_n^{e}\left[q_k(x_1^{-1},\cdots,x_n^{-1})\right]^2=\left[x_1^{\frac{e}{2}} \cdots x_n^{\frac{e}{2}}q_k(x_1^{-1},\cdots,x_n^{-1})\right]^2$$
is a square. Thus, {${\rm SOSrank}(\mathcal{A}) \le n$ in this case}.

Suppose that  Case {\rm (ii)} holds, i.e.,  $n = 3$ and $m=4$. Then by Hilbert's theorem \cite{H88}, ${\rm SOSrank} (\mathcal{A}) \le 3=n$.

Suppose that Case {\rm (iii)} holds, i.e.,  $n=3$ and $m \ge 3e-4$. In the case of $m=en-2=3e-2$ and $m=en=3e$, using similar argument as in the Case {\rm (i)}, we see that the conclusion follows. The only remaining case is when $m=3e-4$. In this case, as $\mathcal{A}$ is a BD-tensor with order $m$, dimension $3$ and exponent $e$ and $\mathcal{A}$ has SOS decomposition, we have
\[
\mathcal{A}{\bf x}^m=x_1^{e}x_2^ex_3^e G(x_1^{-1},x_2^{-1},x_3^{-1}),
\]
where $G$ is a positive semi-definite form and is of 3 dimension and degree $4$. It then from Hilbert's theorem \cite{H88} that $G(t_1,t_2,t_3)$ can be expressed as at most the sum of 3 squares of $3$-dimensional quadratic forms. Thus, using similar line of argument as in Case {\rm (i)} and noting that $e \ge 4$ (as $m=3e-4$ and $m \ge 4$), we have ${\rm SOSrank}(\mathcal{A}) \le n=3$.

Combining these three cases, we see that ${\rm SOSrank}(\mathcal{A}) \le n$, and ${\rm SOSrank}(\mathcal{A})=1$ if $m=en$. In particular, we have
${\rm \mbox {SOS-width}}({\rm BD}_{m,n}^{e} \cap {\rm SOS}_{m,n}) \le n$, and ${\rm \mbox {SOS-width}}({\rm BD}_{m,n}^{e} \cap {\rm SOS}_{m,n}) =1$ if $m=en$.
To see the conclusion, we consider the homogeneous polynomial $$
f_0({\bf x})=\left\{\begin{array}{ll}
x_1^e\cdots x_n^e (\sum_{i=1}^n x_i^{-2}) & \mbox{ if } n \ge 3 \mbox{ and } m = en-2 \\
x_1^2x_2^2+x_2^2x_3^2+x_3^2x_1^2 & \mbox{ if } n=3 \mbox{ and } m=4  \\
x_1^ex_2^ex_3^e(x_1^{-2}x_2^{-2}+x_2^{-2}x_3^{-2}+x_3^{-2}x_1^{-2}) & \mbox{ if } n=3 \mbox{ and } m = 3e-4  \\
\end{array}
\right.$$
and its associated BD-tensor $\mathcal{A}_0$ such that $f_0({\bf x})=\mathcal{A}_0{\bf x}^m$. It can be directly verified that
$${\rm SOSrank}(\mathcal{A}_0)=\left\{\begin{array}{ll}
n & \mbox{ if } n \ge 3 \mbox{ and } m = en-2, \\
3 & \mbox{ if } n=3 \mbox{ and } m=4,  \\
3 & \mbox{ if } n=3 \mbox{ and } m = 3e-4.  \\
\end{array}
\right.$$
For example, in the case $n \ge 3$ and $m = en-2$, to see ${\rm SOSrank}(\mathcal{A}_0)=n$, we only need to show ${\rm SOSrank}(\mathcal{A}_0) \ge n$. Suppose on the contrary that ${\rm SOSrank}(\mathcal{A}_0) \le n-1$. Then, there exists $r \le n-1$ and homogeneous polynomial $f_i$ with degree $m/2=\frac{e}{2}n-1$ such that
\[
x_1^e\cdots x_n^e \left(\sum_{i=1}^n x_i^{-2}\right)=\sum_{i=1}^r f_i({\bf x})^2.
\]
This implies that for each ${\bf x}=(x_1,\cdots,x_n)$ with $x_i \neq 0$, $i=1,\cdots,n$
\[
\sum_{i=1}^n x_i^{-2} = \sum_{i=1}^r \left[\frac{f_i({\bf x})}{x_1^{\frac{e}{2}} \cdots x_n^{\frac{e}{2}}}\right]^2
\]
Letting $t_i=x_i^{-1}$, by continuity, we see that the quadratic form $\sum_{i=1}^n t_i^2$ can be written as a sum of at most $r$ many squares of rational functions in $(t_1,,\cdots,t_n)$. Then, the Cassels-Pfister's Theorem \cite[Theorem 17.3]{book0} (see also \cite[Corollary 17.6]{book0}), implies that the quadratic form $\sum_{i=1}^n t_i^2$ can be written as a sum of at most $r$ many sums of squares of polynomial functions in $(t_1,,\cdots,t_n)$, which is impossible.

In the case $n =3$ and $m = 4$,  we only need to show ${\rm SOSrank}(\mathcal{A}_0) \ge 3$. Suppose on the contrary that ${\rm SOSrank}(\mathcal{A}_0) \le 2$.  Then, there exist $a_i,b_i,c_i,d_i,e_i,f_i \in \mathbb{R}$, $i=1,2$, such that \begin{eqnarray*}
x_1^2x_2^2+x_2^2x_3^2+x_3^2x_1^2&= & (a_1 x_1^2+b_1 x_2^2 + c_1x_3^2 +d_1 x_1x_2+e_1 x_1x_3+f_1x_2x_3)^2 \\
& & +(a_2 x_1^2+b_2 x_2^2 + c_2x_3^2 +d_2 x_1x_2+e_2 x_1x_3+f_2x_2x_3)^2.
\end{eqnarray*}
Comparing with the coefficients gives us that $a_1=a_2=b_1=b_2=c_1=c_2=0$ and
\[
\left\{\begin{array}{l}
d_1^2+d_2^2=1 \\
e_1^2+e_2^2=1 \\
f_1^2+f_2^2=1 \\
d_1e_1+d_2e_2=0 \\
d_1f_1+d_2f_2=0 \\
e_1f_1+e_2f_2=0.
\end{array} \right.
\]
From the last three equations, we see that one of $d_1,d_2,e_1,e_2,f_1,f_2$ must be zero. Let us assume say $d_1=0$. Then, the first equation shows $d_2=\pm 1$ and hence, $e_2=0$ (by the fourth equation). This implies that $e_1= \pm 1$ and $f_2=0$. Again, we have $f_1=\pm 1$ and hence
\[
e_1f_1+e_2f_2=(\pm 1)(\pm 1)+0=\pm 1 \neq 0.
\]
This leads to a contradiction.

For the last case, suppose again by contradiction that ${\rm SOSrank}(\mathcal{A}_0) \le 2$. Then, there exist two homogeneous polynomial $f_i$ with degree $m/2=\frac{3e}{2}-2$ such that
\[
x_1^ex_2^ex_3^e(x_1^{-2}x_2^{-2}+x_2^{-2}x_3^{-2}+ x_3^{-2}x_1^{-2})=\sum_{i=1}^2 f_i({\bf x})^2.
\]
This implies that for each $x=(x_1,\cdots,x_n)$ with $x_i \neq 0$, $i=1,\cdots,n$
\[
x_1^{-2}x_2^{-2}+x_2^{-2}x_3^{-2}+x_3^{-2}x_1^{-2} = \sum_{i=1}^2 \left[\frac{f_i({\bf x})}{x_1^{\frac{e}{2}} \cdots x_3^{\frac{e}{2}}}\right]^2
\]
Letting $t_i=x_i^{-1}$, using a similar line argument in the case $m=en-2$, we see that the polynomial $t_1^2t_2^2+t_2^2t_3^2+t_3^2t_1^2$ can be written
as sums of $2$ squares of polynomials in $(t_1,t_2,t_3)$. This is impossible by the preceding case. Therefore, the conclusion follows. \qed
\end{proof}

Below, let us mention that calculating the exact SOS-rank of SOS tensor decomposition is not a trivial task even for the identity tensor, and this relates to
some open question in algebraic geometry in the literature.  To explain this, we recall that
the identity tensor $\mathcal{I}$ with order $m$ and dimension $n$ is given by $\mathcal{I}_{i_1\cdots i_m}=1$ if $i_1=\cdots=i_m$ and $\mathcal{I}_{i_1\cdots i_m}=0$ otherwise. The identity tensor $\mathcal{I}$ induces the polynomial $f_{\mathcal{I}}({\bf x})=\mathcal{I}{\bf x}^m=x_1^m+\cdots +x_n^m$. It is clear that, $\mathcal{I}$ has an SOS tensor decomposition when $m$ is even and the corresponding SOS-rank of $\mathcal{I}$ is less than or equal to $n$. It was conjectured by Reznick \cite{Rez0} that $f_{\mathcal{I}}({\bf x})$ cannot be written as sums of $(n-1)$ many squares, that is, ${\rm SOSrank}(\mathcal{I})=n$. The positive answer for this conjecture in the special case
of $m=n=4$ was provided in \cite{JPAA1,JPAA2}. On the other hand, the answer for this conjecture in the general case is still open to the best of our knowledge.  Moreover, this conjecture relates to another conjecture of Reznick \cite{Rez0} in the same paper where he showed that the polynomial $f_R({\bf x})=x_1^n+\cdots +x_n^n- n x_1 \cdots x_n$ can be written as sums of $(n-1)$ many squares whenever $n=2^k$ for some $k \in \mathbb{N}$, and he
 conjectured that the estimate of the numbers of squares is sharp. Indeed, he also showed that this conjecture is true whenever the previous conjecture of ``$f_{\mathcal{I}}({\bf x})$ cannot be written as sums of $(n-1)$ many squares" is true.

\setcounter{equation}{0}
\section{Applications}
In this section, we provide some applications for the SOS tensor decomposition of the structure tensors such as finding the minimum $H$-eigenvalue of an
even order extended $Z$-tensor and testing the positive definiteness of a multivariate form. We also provide some numerical examples/experiments to support the theoretical findings.
Throughout this section, all numerical experiments are performed on a desktop, with  3.47 GHz quad-core Intel E5620 Xeon 64-bit CPUs and 4 GB RAM,
equipped with  Matlab 2015.

\subsection{Finding the minimum $H$-eigenvalue of an even order symmetric extended $Z$-tensor}
Finding the minimum eigenvalue of a tensor is an important topic in tensor computation and multilinear algebra, and has found numerous applications including automatic control and image processing \cite{Qi05}.
Recently, it was shown that the minimum $H$-eigenvalue of an even order symmetric $Z$-tensor \cite{HLQS,Hu14} can be found by solving a sums-of-squares optimization problem, which can be
equivalently reformulated as a semi-definite programming problem, and so, can be solved efficiently. In \cite{HLQS}, some upper and lower estimates for the minimum $H$-eigenvalue of general symmetric tensors with even order are
provided via sums-of-squares programming problems. Examples show that the estimate can be sharp in some cases.

On the other hand, it was unknown in \cite{Hu14,HLQS} that whether similar results can continue to hold for some classes of symmetric tensors which are not $Z$-tensors, that is, for symmetric tensors
with possible positive entries on the off-diagonal elements. In this section, as applications of the derived SOS decomposition of structured tensors, we show that the class of even order symmetric extended $Z$-tensor serves as one such class. To present the conclusion, the following Lemma plays an important role in our later analysis.
\bl\label{lema51}$^{\cite{Qi05}}$ \label{lemma:5.1} Let $\mathcal{A}$ be a symmetric tensor with even order $m$ and dimension $n$. Denote the
minimum $H$-eigenvalue of $\mathcal{A}$  by
$\lambda_{min}(\mathcal{A})$. Then, we have
\begin{equation}\label{e51}\lambda_{min}(\mathcal{A})=\min_{{\bf x\neq {\bf 0}}}\frac{\mathcal{A}{\bf x}^m}{\|{\bf x}\|_m^m}=\min_{\|{\bf x}\|_m=1}\mathcal{A}{\bf x}^m,
\end{equation}
where $\|{\bf x}\|_m=(\sum_{i=1}^n|x_i|^m)^{\frac{1}{m}}$.
\el

\begin{theorem}{\bf (Finding the minimum $H$-eigenvalue of an even order symmetric extended $Z$-tensor)} \label{th:5.1}
Let $m$ be an even number. Let $\mathcal{A}$ be a symmetric extended $Z$-tensor with order $m$ and dimension $n$. Then, we have
\begin{eqnarray*}
\lambda_{{\rm min}}(\mathcal{A})&=& \max_{\mu, r \in \mathbb{R}}\{\mu: f_{\mathcal{A}}({\bf x})-r (\|{\bf x}\|_m^m-1)-\mu \in \Sigma^2_m[{\bf x}]\},
\end{eqnarray*}
where $f_{\mathcal{A}}({\bf x})=\mathcal{A}{\bf x}^m$ and $\Sigma_m^2[{\bf x}]$ is the set of all SOS polynomials with degree at most $m$.
 \end{theorem}
\proof
Consider the following problem
  \[
(P) \ \ \  \min \{\mathcal{A}{\bf x}^m: \|{\bf x}\|_m^m=1\}
\]
  and denote its global minimizer by ${\bf a}=(a_1,\cdots,a_n)^T \in \mathbb{R}^n$.  Clearly, $\sum_{i=1}^na_i^m=1$. Then, $\lambda_{\rm min}(\mathcal{A})=f_{\mathcal{A}}({\bf a})=\mathcal{A}{\bf a}^m$.
It follows that for all ${\bf x} \in \mathbb{R}^n\backslash\{0\}$
\begin{eqnarray*}
f_{\mathcal{A}}({\bf x})-\lambda_{\rm min}(\mathcal{A}) \sum_{i=1}^nx_i^m & = & f_{\mathcal{A}}({\bf x})-f_{\mathcal{A}}({\bf a})\sum_{i=1}^nx_i^m \\
&= & \sum_{i=1}^nx_i^m\bigg(f_{\mathcal{A}}(\frac{{\bf x}}{(\sum_{i=1}^nx_i^m)^{\frac{1}{m}}})-f_{\mathcal{A}}({\bf a})\bigg) \ge 0,
\end{eqnarray*}
where the last inequality holds as $m$ is even and $\overline{{\bf x}}=\frac{{\bf x}}{(\sum_{i=1}^nx_i^m)^{\frac{1}{m}}}$ belongs to the feasible set of (P).
 This shows that
$g({\bf x}):=f_{\mathcal{A}}({\bf x})-\lambda_{\rm min}(\mathcal{A}) \sum_{i=1}^nx_i^m$ is a homogeneous polynomial which always take nonnegative values. As $\mathcal{A}$ is an extended $Z$-tensor,
 there exist $s \in \mathbb{N}$ and index sets $\Gamma_l \subseteq \{1,\cdots,n\}$, $l=1,\cdots,s$ with  $\bigcup_{l=1}^s \Gamma_l=\{1,\cdots,n\}$ and $\Gamma_{l_1} \cap \Gamma_{l_2} =\emptyset$  such that
for all ${\bf x} \in \mathbb{R}^n$
\begin{equation}\label{eq:96}
f_{\mathcal{A}}({\bf x})=\sum_{i=1}^n f_{m,i} x_i^{m}+\sum_{l=1}^s \sum_{\alpha_l \in \Omega_l}f_{\alpha_l} {\bf x}^{\alpha_l}
\end{equation}
such that, for each $l=1,\cdots,s$, either one of the following two condition holds:
(1) $f_{\alpha_l}=0$ for all but one $\alpha_l \in \Omega_l$; (2) $f_{\alpha_l} \le 0$ for all $\alpha_l \in \Omega_l$.
 Thus,
\[
g({\bf x})=\sum_{i=1}^n (f_{m,i}-\lambda_{\rm min}(\mathcal{A})) x_i^{m}+\sum_{l=1}^s \sum_{\alpha_l \in \Omega_l} f_{\alpha_l} {\bf x}^{\alpha_l},
\]
is an extended $Z$-polynomial which always takes nonnegative values. Let $\mathcal{B}$ be a symmetric tensor such that $g({\bf x})=\mathcal{B}{\bf x}^m$. Then,
$\mathcal{B}$ is a positive semi-definite extended $Z$-tensor and so is SOS by Theorem \ref{th:5.2}. Thus, $g({\bf x})$ is an SOS polynomial with degree $m$. Note that $g({\bf x})=
f_{\mathcal{A}}({\bf x})-\lambda_{\rm min}(\mathcal{A}) \sum_{i=1}^nx_i^m=
f_{\mathcal{A}}({\bf x})-\lambda_{\rm min}(\mathcal{A}) (\sum_{i=1}^nx_i^m-1)-\lambda_{\rm min}(\mathcal{A})$. This shows that
\[
\lambda_{\rm min}(\mathcal{A}) \le \max_{\mu, r \in \mathbb{R}}\{\mu: f_{\mathcal{A}}({\bf x})-r (\|{\bf x}\|_m^m-1)-\mu \in \Sigma^2_m[{\bf x}]\}.
\]
To see the reverse inequality, take any $(\mu,r)$ with $f_{\mathcal{A}}({\bf x})-r (\|{\bf x}\|_m^m-1)-\mu \in \Sigma^2_m[{\bf x}]$. Then, for all ${\bf x} \in \mathbb{R}^n$,
\[
f_{\mathcal{A}}({\bf x})-r (\|{\bf x}\|_m^m-1)-\mu \ge 0.
\]
This shows that $r \ge \mu$ and $f_{\mathcal{A}}({\bf x}) \ge r \|{\bf x}\|_m^m$ for all ${\bf x} \in \mathbb{R}^n$. This shows that $\lambda_{\rm min}(\mathcal{A}) \ge r \ge \mu$, and so, the conclusion follows. \qed

\begin{Remark}\label{remark:5.1}
Let $\mathcal{A}$ be an extended $Z$-tensor. As in (\ref{eq:96}), its associated polynomial $f_{\mathcal{A}}$ can be written as $f_{\mathcal{A}}({\bf x})=\sum_{i=1}^n f_{m,i} x_i^{m}+\sum_{l=1}^s \sum_{\alpha_l \in \Omega_l}f_{\alpha_l} {\bf x}^{\alpha_l}$
. Then, Remark \ref{remark:3.3} implies that
\begin{eqnarray*}
\lambda_{{\rm min}}(\mathcal{A}) & = & \max_{\mu, r \in \mathbb{R}}\{\mu: f_{\mathcal{A}}({\bf x})-r (\|{\bf x}\|_m^m-1)-\mu \in \Sigma^2_m[{\bf x}]\} \\
& = & \max_{\mu, r \in \mathbb{R}}\{\mu: f_{\mathcal{A}}({\bf x})-r \|{\bf x}\|_m^m \in \Sigma^2_m[{\bf x}], r-\mu \ge 0\} \\
& = & \max_{\mu, r \in \mathbb{R}}\{\mu: \sum_{i \in \Gamma_l} f_{m,i} x_i^{m}+ \sum_{\alpha_l \in \Omega_l}f_{\alpha_l} {\bf x}^{\alpha_l}-r \|{\bf x}^{(l)}\|_m^m \in \Sigma^2_m[{\bf x}^{(l)}], l=1,\cdots,s \\
& & \ \ \ \ \ \ \ \ \ \ \ \ \ r-\mu \ge 0\}
\end{eqnarray*}
where, for each $l=1,\cdots,s$, ${\bf x}^{(l)}=(x_i)_{i \in \Gamma_l}$ and $\Sigma^2_m[{\bf x}^{(l)}]$ is the set of all
SOS polynomials in ${\bf x}^{(l)}$.
\end{Remark}

As explained in \cite{Hu14,HLQS}, the sums-of-squares problem $$\max_{\mu, r \in \mathbb{R}}\{\mu: f_{\mathcal{A}}({\bf x})-r (\|{\bf x}\|_m^m-1)-\mu \in \Sigma^2_m[{\bf x}]\}$$ can be equivalently rewritten as a semi-definite programming problem (SDP), and so, can be solved efficiently. Indeed, this conversion can be done by using the commonly used
Matlab Toolbox YALMIP \cite{Yalmip1,Yalmip2}.
On the other hand, the size of the equivalent SDP problem of the relaxation problem increase dramatically when the  dimension/order of the tensor increases.
For example, as illustrate in Table 1, for a $4$th-order $50$-dimensional tensor, the equivalent SDP problem has $1326$ variables and $316251$ constraints.  Fortunately, a robust SDP software (SDPNAL \cite{zhao}) has been established recently which enables us to solve large-scale SDP (dimension up to 5000 and number of constraint of the SDP up to 1 million). This
enables us to find the minimum $H$-eigenvalue for medium-size tensor. Later on, we will explain how to use SDPNAL together with the observation in Remark \ref{remark:5.1}
to find the minimum $H$-eigenvalue for large-size tensor.

We first illustrate how to compute the minimum $H$-eigenvalue of an extended $Z$-tensor $\mathcal{A}$ using the above sums-of-squares problem via Matlab Toolbox YALMIP \cite{Yalmip1,Yalmip2} via
two small-size problems. We will show the performance of the method for various larger-size problem later.
\begin{example}
Consider the symmetric tensor $\mathcal{A}$ with order $6$ and dimension $4$ where
\[
\mathcal{A}_{111111}=\mathcal{A}_{222222}=\mathcal{A}_{333333}=\mathcal{A}_{444444}=1,
\]
\[
\mathcal{A}_{i_1 \cdots i_6}=\frac{1}{5}, \mbox{ for all } (i_1,\cdots,i_6)=\sigma(1,1,1,2,2,2),
\]
\[
\mathcal{A}_{i_1 \cdots i_6}=\frac{2}{5}, \mbox{ for all } (i_1,\cdots,i_6)=\sigma(3,3,4,4,4,4),
\]
and $\mathcal{A}_{i_1 \cdots i_6}=0$ otherwise. Here $\sigma(i_1,\cdots,i_6)$ denotes all the possible permutation of $(i_1,\cdots,i_6)$. The associated polynomial
$$f_{\mathcal{A}}({\bf x})=\mathcal{A}{\bf x}^m=x_1^6+x_2^6+x_3^6+x_4^6+4x_1^3x_2^3+6x_3^2x_4^4$$ is an extended  $Z$-polynomial. So, $\mathcal{A}$ is an extended $Z$- tensor. It can be easily verified that $\mathcal{A}$ is not a $Z$-tensor.

To compute its minimum $H$-eigenvalue, we note that the corresponding sums-of-squares optimization problem reads
  \[
  \max_{\mu, r \in \mathbb{R}}\{\mu: f_{\mathcal{A}}({\bf x})-r (\|{\bf x}\|_6^6-1)-\mu \in \Sigma^2_6[{\bf x}]\}.
  \]
Convert this sums-of-squares optimization problem into a semi-definite programming problem using the Matlab Toolbox {\rm YALMIP} \cite{Yalmip1,Yalmip2}, and solve it by using the
SDP software {\rm SDPNAL} we obtain that $\lambda_{\rm min}(\mathcal{A})=-1$.
The simple code using {\rm YALMIP} is appended as follows:
\begin{verbatim}
sdpsettings('solver','sdpnal')
sdpvar x1 x2 x3 x4 r mu
f = x1^6+x2^6+x3^6+x4^6+4*x1^3*x2^3+6*x3^2*x4^4;
g = [(x1^6+x2^6+x3^6+x4^6)-1];
F = [sos(f-mu-r*g)];
solvesos(F,-mu,[],[r;mu])
\end{verbatim}

Moreover, note from the geometric mean inequality that
$|x_1^3x_2^3|=(x_1^6)^{\frac{1}{2}} (x_2^6)^{\frac{1}{2}}\le \frac{1}{2} x_1^6+\frac{1}{2} x_1^6$.
It follows that
\[
f_{\mathcal{A}}({\bf x})+\|{\bf x}\|_6^6=2x_1^6+2x_2^6+2x_3^6+2x_4^6+4x_1^3x_2^3+6x_3^2x_4^4 \ge 0 \mbox{ for all } {\bf x} \in \mathbb{R}^n.
\]
On the other hand, consider $\bar {\bf x}=(\sqrt[6]{\frac{1}{2}},-\sqrt[6]{\frac{1}{2}},0,0)$. We see that $f_{\mathcal{A}}(\bar {\bf x})+\|\bar {\bf x}\|_6^6=0$. This shows that $\lambda_{\rm min}(\mathcal{A})=\min\{f_{\mathcal{A}}({\bf x}): \|{\bf x}\|_6=1\}=-1$. This verifies the correctness of our computed minimum $H$-eigenvalue.

\end{example}

\begin{example}
Let $\alpha,\beta \in \mathbb{R}$ and consider the symmetric tensor $\mathcal{A}$ with order $6$ and dimension $4$ where
\[
\mathcal{A}_{111111}=\mathcal{A}_{222222}=\mathcal{A}_{333333}=\mathcal{A}_{444444}=1,
\]
\[
\mathcal{A}_{i_1 \cdots i_6}=\alpha, \mbox{ for all } (i_1,\cdots,i_6)=\sigma(1,1,1,2,2,2),
\]
\[
\mathcal{A}_{i_1 \cdots i_6}=\beta, \mbox{ for all } (i_1,\cdots,i_6)=\sigma(3,3,3,4,4,4),
\]
and $\mathcal{A}_{i_1 \cdots i_6}=0$ otherwise. Here $\sigma(i_1,\cdots,i_6)$ denotes all the possible permutation of $(i_1,\cdots,i_6)$. The associated polynomial
$$f_{\mathcal{A}}({\bf x})=\mathcal{A}{\bf x}^m=x_1^6+x_2^6+x_3^6+x_4^6+20 \alpha \,x_1^3x_2^3+20 \beta\, x_3^3x_4^3$$ is an extended  $Z$-polynomial. So, $\mathcal{A}$ is an extended $Z$- tensor. It can be easily verified that if either $\alpha>0$ or $\beta>0$, then $\mathcal{A}$ is not a $Z$-tensor.

To compute its minimum $H$-eigenvalue, we randomly generate 100 instance of $(\alpha,\beta) \in [-5,5] \times [-5,5]$. For each $(\alpha,\beta)$, we  convert  the corresponding sums-of-squares optimization problem
  \[
  \max_{\mu, r \in \mathbb{R}}\{\mu: f_{\mathcal{A}}({\bf x})-r (\|{\bf x}\|_6^6-1)-\mu \in \Sigma^2_6[{\bf x}]\}
  \]
into a semi-definite programming problem using the Matlab Toolbox {\rm YALMIP} \cite{Yalmip1,Yalmip2}, and solve it by using the
SDP software {\rm SDPNAL}.  We then compare the computed minimum $H$-eigenvalue with the true minimum $H$-eigenvalue of $\mathcal{A}$. Indeed,  similar to the preceding example, we can verify that $\lambda_{\rm min}(\mathcal{A})=m(\alpha,\beta)$ where
\[
m(\alpha,\beta):=\left\{
\begin{array}{cll}
1-10|\alpha| & \mbox{ if } &  |\alpha| \ge  |\beta|,\\

1-10 |\beta| & \mbox{ if } & |\alpha| <    |\beta|.
\end{array}
 \right.
\]


For all the $100$ generated $(\alpha,\beta)$, the maximum difference of the computed $H$-minimum eigenvalue and the true $H$-minimum eigenvalue is
$6.2039e-05$.
\end{example}

\subsection*{Medium-size examples}
We now consider a few medium-size examples which involves symmetric extended $Z$-tensor with order up to $30$ or dimension up to $60$ .
\begin{example}
Let $m=10k$ with $k \in \mathbb{N}$. Consider the symmetric tensor $\mathcal{A}$ with order $m$ and dimension $4$ where
\[
\mathcal{A}_{1 \cdots 1}=\mathcal{A}_{2 \cdots 2}=\mathcal{A}_{3 \cdots 3}=\mathcal{A}_{4 \cdots 4}=1,
\]
\[
\mathcal{A}_{i_1 \cdots i_{m}}=\alpha, \mbox{ for all } (i_1,\cdots,i_{m})=\sigma(\underbrace{1,\cdots,1}_{m/2},\underbrace{2,\cdots,2}_{m/2}),
\]
\[
\mathcal{A}_{i_1 \cdots i_{m}}= \beta, \mbox{ for all } {  (i_1,\cdots,i_{m})=\sigma(\underbrace{3,\cdots,3}_{m/5},\underbrace{4,\cdots,4}_{4m/5}),}
\]
\[
\mathcal{A}_{i_1 \cdots i_{m}}= \beta, \mbox{ for all } {  (i_1,\cdots,i_{m})=\sigma(\underbrace{3,\cdots,3}_{4m/5},\underbrace{4,\cdots,4}_{m/5}),}
\]
with $\alpha=2 {m \choose m/2} ^{-1}$ and $\beta=- {m \choose m/5} ^{-1}$,
and $\mathcal{A}_{i_1 \cdots i_m}=0$ otherwise. Here $\sigma(i_1,\cdots,i_m)$ denotes all the possible permutation of $(i_1,\cdots,i_m)$. The associated polynomial
$${  f_{\mathcal{A}}({\bf x})=\mathcal{A}{\bf x}^m=x_1^{m}+x_2^{m}+x_3^{m}+x_4^{m}+ 2x_1^{\frac{m}{2}}x_2^{\frac{m}{2}}-x_3^\frac{m}{5}x_4^{\frac{4m}{5}}-x_3^{\frac{4m}{5}}x_4^{\frac{m}{5}},}$$
is an extended  $Z$-polynomial. So, $\mathcal{A}$ is an extended $Z$-tensor. It can be easily verified that $\mathcal{A}$ is not a $Z$-tensor.
Moreover, using geometric mean inequality, we can  directly verify that the true minimum $H$-eigenvalue is $0$.

We compute the minimum
$H$-eigenvalue by solving the corresponding sums-of-squares problem for the case $m=20, 30$, and compare with the true minimum $H$-eigenvalue. The results are
summarized in Table 1.
%
%
\end{example}

\begin{example}
Let $n=4k$ with $k \in \mathbb{N}$. Consider the symmetric tensor $\mathcal{A}$ with order $4$ and dimension $n$ where
\[
\mathcal{A}_{1111}=\mathcal{A}_{2222}=\cdots=\mathcal{A}_{nnnn}=n,
\]
\[
\mathcal{A}_{i_1i_2i_3i_4}=\frac{1}{6}, \mbox{ for all } (i_1,i_2,,i_3,i_4)=\sigma(4i-3,4i-2,4i-1,4i), i=1,\cdots,\frac{n}{4},
\]
and $\mathcal{A}_{i_1i_2i_3i_4}=0$ otherwise. Here $\sigma(i_1,\cdots,i_4)$ denotes all the possible permutation of $(i_1,\cdots,i_4)$. The associated polynomial
$$f_{\mathcal{A}}({\bf x})=\mathcal{A}{\bf x}^m= n(x_1^{4}+\cdots+x_n^4) + 4\sum_{i=1}^{n/4}x_{4i-3}\, x_{4i-2}\, x_{4i-1}\, x_{4i}$$  is an extended  $Z$-polynomial. So, $\mathcal{A}$ is an extended $Z$-tensor. It can be easily verified that $\mathcal{A}$ is not a $Z$-tensor.
Moreover, using geometric mean inequality, we can directly verify that the true minimum $H$-eigenvalue is $n-1$.

We compute the minimum
$H$-eigenvalue by solving the corresponding  sums-of-squares problem for the case $n=20,40,50,60$, and compare with the true minimum $H$-eigenvalue. The results are
summarized in Table 1.

%
%
\end{example}

The following table summarizes the numerical results of Example 5.3 and Example 5.4 where we compute the minimum $H$-eigenvalue by first converting the corresponding sums-of-squares problem to
an SDP problem using YALMIP and solving this SDP problem using SDPNAL.  In particular,  the data of the following table are explained as follows. \begin{itemize}
\item $m$: the order of the symmetric tensor,
\item $n$: the dimension of the symmetric tensor,
 \item $NV$: the number of variables of the equivalent SDP problem,
 \item $NC$: the number of constraints in the equivalent SDP problem,
 \item Computed eigenvalue: the calculated minimum $H$-eigenvalue,
 \item True eigenvalue: the true minimum $H$-eigenvalue
 \item Time (YALMIP): the CPU-time for converting the sums-of-squares problem to SDP (measured in seconds).
 \item Time (SDPNAL): the CPU-time for solving SDP via SDPNAL (measured in seconds).
 \end{itemize}

\begin{table}[h!]

\begin{center}

\caption{Test results for medium size tensors}
\begin{tabular}{|c||c|c|c|c|c|c|c|c|}\hline
Problem & m & n &   NV & NC & Computed  & True  & Time (YALMIP) & Time (SDPNAL) \\
& & & & & eigenvalue & eigenvalue & & \\ \hline
Example 5.3 & 20 & 4 & 1001 & 1001 & -1.7634e-09 & 0 & 11.9487 & 0.5700 \\   \hline
Example 5.3 & 30 & 4 & 3876 & 6936 & 1.1382e-12 & 0 &  198.8141 & 8.2700 \\   \hline
Example 5.4 & 4 & 20 & 231 & 10626 & 19.0000 & 19 & 4.6951 & 0.4763 \\   \hline
Example 5.4 & 4 & 40 & 861 &  135751 & 39.0000  & 39 & 440.8231 & 1.7727 \\   \hline
Example 5.4 & 4 & 50  &  1326 &  316251 & 49.0000 & 49 & 2365.9043& 5.1109 \\   \hline
Example 5.4 & 4 & 60 & 1891 & 635376 & 59.0000 & 59 & 9322.0631 & 50.2934 \\ \hline
\end{tabular}
\end{center}
\end{table}
We observe that, for all the above numerical examples, the minimum $H$-eigenvalues can be found successfully
for medium-size tensors.

\subsection{Large size examples}
Finally, we illustrate with an example that using SDPNAL together with the observation in Remark \ref{remark:5.1} enables us to solve
some large size tensors (dimension up to 2000).

As one can observed in Table 1, most of the time are occupied in YALMIP in converting the sums-of-squares problem into an SDP problem. This process involves matching up
the coefficients of all the involved ${m+n-1 \choose m}$ many monomials, and so, can be time-consuming. On the other hand, by using
the sums-of-squares problem discussed in Remark \ref{remark:5.1} and letting $k=\max_{1 \le l \le s}|\Gamma_l|$, the corresponding process only involves  $s {m+k-1 \choose m}$ many monomials which is much smaller than ${m+n-1 \choose m}$ when
$s$ is large and $k$ is small. For example, as in Example 5.4, we can set $s=n/4$, $k=4$ and $m=4$, and so, $s {m+k-1 \choose m}$ is of the order $n$; while ${m+n-1 \choose m}={n+3 \choose 4}$ which
is of the order $n^4$.

The following table summarizes the numerical results of  Example 5.4 with dimension from 500 to 2000, where we compute the minimum $H$-eigenvalue by first converting the corresponding sums-of-squares problem
discussed in Remark \ref{remark:5.1} to
an SDP problem using YALMIP and solving this SDP problem using SDPNAL. We observe that, for all the instances, the minimum $H$-eigenvalues can be found successfully.
The meaning of the data are the same as in Table 1.
\begin{table}[h!]
\begin{center}
\caption{Test results for large size tensors}
\begin{tabular}{|c||c|c|c|c|c|c|c|c|}\hline
Problem & m & n &   NV & NC & Computed  & True  & Time (YALMIP) & Time (SDPNAL) \\
& & & & & eigenvalue & eigenvalue & & \\ \hline
Example 5.4 & 4 & 500 & 1250 & 1375 & 499.0000 & 499 & 4.6299 & 6.8295 \\   \hline
Example 5.4 & 4 & 1000 & 2500 &  2750 & 999.0000  & 999 & 8.8298 & 66.5566 \\   \hline
Example 5.4 & 4 & 2000  & 5000 & 5500  & 1999.0000 & 1999 & 20.9729 & 563.6903 \\   \hline
\end{tabular}
\end{center}
\end{table}


\newpage
\subsection{Testing positive definiteness of a multivariate form}
For a multivariate form $\mathcal{A}{\bf x}^m$, we say it is positive definite if $\mathcal{A}{\bf x}^m>0$ for all ${\bf x} \neq {\bf 0}$.
Testing positive definiteness of a multivariate form $\mathcal{A}{\bf x}^m$ is an important
problem in the stability study of nonlinear autonomous systems via Lyapunov's direct
method in automatic control \cite{Qi05}. Researchers in automatic control have studied the
conditions of such positive definiteness intensively. However, for $n \ge 3$
and $m \ge 4$, this is, in general, a hard problem in mathematics. Recently, some efficient methods based on eigenvalues
of tensors were proposed to solve the problem in the case where $m=4$ \cite{Ni}.

In this part, we show that testing positive definiteness of a multivariate form $\mathcal{A}{\bf x}^m$ where $\mathcal{A}$ is an extended $Z$-tensor
can be  computed by sums-of-squares problem via Theorem \ref{th:5.1}.  Indeed, a direct consequence of Theorem \ref{th:5.1} and Lemma \ref{lemma:5.1} give us the following useful test:
\begin{corollary}
Let $\mathcal{A}$ be an extended $Z$-tensor. Then, the associated multivariate form $\mathcal{A}{\bf x}^m$ is positive definite if and only if
\[
\max_{\mu, r \in \mathbb{R}}\{\mu: f_{\mathcal{A}}({\bf x})-r (\|{\bf x}\|_m^m-1)-\mu \in \Sigma^2_m[{\bf x}]\}>0,
\]
where $f_{\mathcal{A}}({\bf x})=\mathcal{A}{\bf x}^m$ and $\Sigma_m^2[{\bf x}]$ is the set of all SOS polynomials with degree at most $m$.
\end{corollary}

We now use the above corollary to test the positive definiteness of extended $Z$-tensors. To do this, we first generate 100 extended $Z$-tensors as numerical examples. These extended $Z$-tensors are randomly generated by the
following procedure.
\bigskip

{\bf Procedure 1}
\begin{itemize}
\item[(i)] Given $(m, n,s,k,M)$ with $m$ is an even number and $n=sk$, where $n$ and $m$ are the dimension and the order of the
randomly generated tensor, respectively,  and $M$ is a large positive constant.

\item[(ii)] Randomly generate a random positive integer $L$ and a partition of the index set $\{1,\cdots,n\}$, $\{\Gamma_1,\cdots,\Gamma_s\}$, such that $|\Gamma_i|=k$, $i=1,\cdots,s$ and $\Gamma_{i} \cap \Gamma_{i'}=\emptyset$ for all $i \neq i'$.
For each $i=1,\cdots,s-1$, generate a random multi-index $(l_1^i,\cdots,l_m^i)$ with $l_j^i \in \Gamma_i$, $j=1,\cdots,m$ and a random number $\bar{a}_{l_1^i \cdots l_m^i} \in [0,1]$.  { Generate one randomly $m$th-order $k$-dimensional symmetric tensor $\mathcal{B}$,  such that all elements
of $\mathcal{B}$ are in the interval $[0, 1]$.}

\item[(iii)]  We define extended $Z$-tensor $\mathcal{A}=(a_{i_1i_2\cdots i_m})$ such that
\[
a_{i_1 \cdots i_m}= \left\{\begin{array}{cll}
(-1)^L M  & \mbox{ if } & i_1=\cdots=i_m=i \mbox{ for all } i=1,\cdots,n, \\
\bar{a}_{l_1^i \cdots l_m^i} & \mbox{ if } & (i_1,\cdots,i_m)=\sigma(l_1^i,\cdots,l_m^i) \mbox{ with } l_1^i,\cdots,l_m^i \in \Gamma_i, i=1,\cdots,s-1, \\
-\mathcal{B}_{i_1 \cdots i_m} & \mbox{ if } & {  i_1,\cdots,i_m \in \Gamma_s,}\\
0 & \mbox{ othewise. } &
\end{array}
\right.
\]
Here $\sigma(i_1,\cdots,i_m)$ denotes all the possible permutation of $(i_1,\cdots,i_m)$.
\end{itemize}

From the construction of $\mathcal{A}$, it can be verified that $\mathcal{A}$ is an extended $Z$-tensor.
Let $f_{\mathcal{A}}({\bf x})=\mathcal{A}{\bf x}^m$. We then solve the sums-of-squares problem $$\max_{\mu, r \in \mathbb{R}}\{\mu: f_{\mathcal{A}}({\bf x})-r (\|{\bf x}\|_m^m-1)-\mu \in \Sigma^2_m[{\bf x}]\}$$
and use the preceding corollary to determine whether $\mathcal{A}{\bf x}^m$ is a positive definite multivariate form or not. Here, to speed up the algorithm, as we did for the large size
tensors, we first convert the sums-of-squares problem into an SDP by using Remark \ref{remark:5.1} and {\rm YALMIP}. Then, we solve the equivalent
SDP by using the software {\rm SDPNAL}. The correctness can be verified by looking at the randomly generated positive number $L$. Indeed, from the construction, if $L$ is an even number and $M$ is a large positive number, the diagonal elements will strictly dominate the sum of the off-diagonal elements, and so, $\mathcal{A}{\bf x}^m$ is a positive definite multivariate form. On the other hand, if $L$ is an odd number, then the diagonal elements will be negative, and so, $\mathcal{A}{\bf x}^m$ is not a positive definite multivariate form in this case.

The following table summarize the results for the correctness of testing the positive definiteness of a multivariate form generated by an extended $Z$-tensor.
As we can see the results, in our numerical experiment, all the $100$ randomly generated instance has been correctly identified.
\el

\begin{center}
\begin{tabular}{|c|c|c|c|c|c|c|c|}
  \hline
  m & n & s & k & M & PD & NPD & Correctness \\ \hline
  4 & 20 & 4 & 5 & 100 &  48 & 52 & 100\% \\
  4 & 25 & 5 & 5 & 100 &  46 & 54 &100\% \\
  4 & 40 & 4 & 10 & 100 & 52 & 48 & 100\% \\
  4 & 60 & 4 & 15 & 100 &  45 & 55 & 100\% \\
  4 & 100 & 4 & 25 & 100 & 44 & 56 & 100\% \\
    \hline
\end{tabular}
\end{center}










\section{Conclusions and Remarks}

In this paper, we establish SOS tensor decomposition of various even order symmetric structured tensors available in the current literature.
These include positive Cauchy  tensors, weakly diagonally dominated tensors,
$B_0$-tensors, { double $B$-tensors, quasi-double $B_0$-tensors, $MB_0$-tensors,} $H$-tensors, absolute tensors
of positive semi-definite $Z$-tensors and extended $Z$-tensors.
We also examine the SOS-rank of SOS tensor decomposition and the SOS-width for SOS tensor cones. In particular, we provide an explicit sharp
estimate for SOS-rank of tensors with bounded exponent and SOS-width for the tensor cone consisting of all such tensors with bounded exponent that have SOS decomposition.
It is shown that the SOS-rank of SOS tensor decomposition is equal to the
optimal value of a related rank optimization problem over positive semi-definite matrix constraints. Finally, applications for the SOS decomposition
of extended $Z$-tensors are provided and several numerical experiments illustrate the significance.

Below, we raise some open questions which might be interesting for future work:


{\bf Question 1}: Can we evaluate the SOS-rank of symmetric $B_0$-tensors?

{\bf Question 2}: Can we evaluate the SOS-rank of symmetric $Z$-tensors?

{\bf Question 3}: Can we evaluate the SOS-rank of symmetric diagonally dominated tensors?

{ {\bf Question 4}: Can we use the techniques in Section 5 to find the minimum $H$-eigenvalue
of an even order symmetric structured tensors other than the extended $Z$-tensors? }

\bigskip

\noindent
{\bf Acknowledgment}  We are thankful to Prof. Man-Duen Choi, Prof. Changqing Xu, Dr. Xin Liu and Dr. Ziyan Luo for their comments, which improved our paper.

\end{document}